\documentclass{tran-l}

 \newtheorem{theorem}{Theorem}[subsection]
 \newtheorem{cor}[theorem]{Corollary}
 \newtheorem{lemma}[theorem]{Lemma}
 \newtheorem{proposition}[theorem]{Proposition}
 \theoremstyle{definition}
 \newtheorem{definition}[theorem]{Definition}
 \theoremstyle{definition}
 \newtheorem{example}[theorem]{Example}
 \theoremstyle{remark}
 \newtheorem{rem}[theorem]{Remark}
 \numberwithin{equation}{subsection}

\usepackage{latexsym}
\usepackage{amssymb}

\usepackage[all]{xy}

\newcommand{\ben}{\begin{equation}}
\newcommand{\een}{\end{equation}}

\newcommand{\integer}{\ensuremath{{\mathbf Z}}}
\newcommand{\naturals}{\ensuremath{{\mathbf N}}}
\newcommand{\real}{\ensuremath{{\mathbf R}}}
\newcommand{\complex}{\ensuremath{{\mathbf C}}}

\newcommand{\rational}{\ensuremath{{\mathbf Q}}}

\newcommand{\U}[1]{\ensuremath{{\mathrm U( #1 )}}}
\newcommand{\Or}[1]{\ensuremath{{\mathrm{O}( #1 )}}}

\newcommand{\GL}[1]{\ensuremath{{\mathrm {GL}_{ #1 }}}}
\newcommand{\GLC}[1]{\GL{#1}(\complex)}
\newcommand{\SL}[1]{\ensuremath{{\mathrm {SL}_{ #1 }}}}
\newcommand{\SLC}[1]{\SL{#1}(\complex)}
\newcommand{\GLres}{\ensuremath{{\mathrm{GL}}_{\mathrm{res}}(\HH)}}

\newcommand{\Aa}{{\mathcal A}}
\newcommand{\XX}{{\mathcal X}}
\newcommand{\BB}{{\mathcal B}}
\newcommand{\KK}{{\mathcal K}}

\newcommand{\UU}{{\mathcal U}}
\newcommand{\VV}{{\mathcal V}}

\newcommand{\SSS}{{\mathcal S}}
\newcommand{\FF}{{\mathcal F}}
\newcommand{\GG}{{\mathcal G}}
\newcommand{\CC}{\mathcal{C}}

\newcommand{\LL}{\mathcal{L}}
\newcommand{\MM}{\mathcal{M}}
\newcommand{\HH}{\mathcal{H}}
\newcommand{\RR}{\mathcal{R}}
\newcommand{\PP}{\mathcal{P}}
\newcommand{\Proj}{\mathbf{P}}

\newcommand{\Xx}{\mathsf{X}}
\newcommand{\Loop}{\mathsf{L}}

\newcommand{\twoarrows}{\rightrightarrows}

\newcommand{\To}{\longrightarrow}

\newcommand{\timests}{\: {}_{t}  \! \times_{s}}
\newcommand{\fiberprod}[2]{\: {}_{#1}  \! \times_{#2}}
\newcommand{\toparrow}[1]{\stackrel{#1}{\rightarrow}}
\newcommand{\Gpd}{\RR \rightrightarrows \UU}
\newcommand{\Gerbe}{{\texttt{Gb}}}
\newcommand{\Semiproductofin}{{\U{n} \widetilde{\times}
\Proj\U{n}}}
\newcommand{\UK}{{\mathrm{U}}_\KK}
\newcommand{\Semiproducto}{{\UK \widetilde{\times} \Proj
\U{\HH}}}

\newcommand{\Grr}{\ensuremath{{\mathrm{Gr}}}}

\newcommand{\Korb}[1]{\ensuremath{{}^{#1} K_{\mathrm{orb}}}}
\newcommand{\Kgpd}[1]{\ensuremath{{}^{#1} K_{\mathrm{gpd}}}}

\begin{document}

\title[Gerbes over Orbifolds and K-theory]{Gerbes over Orbifolds and  Twisted K-theory}
\author{ Ernesto Lupercio and Bernardo Uribe}

\address{Department of Mathematics, University of Wisconsin at Madison, Madison, WI 53706}
\address{Department of Mathematics, University of Wisconsin at Madison, Madison, WI 53706}
\email{ lupercio@math.wisc.edu \\ uribe@math.wisc.edu}

\maketitle

\begin{abstract}
   In this paper we construct an explicit geometric model for the group of gerbes
   over an orbifold $X$. We show how from its curvature we can obtain its characteristic
   class in $H^3(X)$ via Chern-Weil theory.
   For an arbitrary gerbe $\LL$, a twisting $\Korb{\LL}(X)$ of the orbifold
   $K$-theory of $X$ is constructed, and shown to generalize previous twisting by Rosenberg
   \cite{Rosenberg}, Witten \cite{Witten}, Atiyah-Segal \cite{Atiyah} and Bowknegt et.
   al. \cite{Murray2}
   in the smooth case and by Adem-Ruan \cite{AdemRuan} for discrete torsion
   on an orbifold.
\end{abstract}

\maketitle

\tableofcontents

\section{Introduction}

   An orbifold is a very natural generalization of a manifold. Locally it looks like the
quotient of an open set of a vector space divided by the action of a group, in such a
way that  the stabilizer of the action at every point is a finite group. Many moduli
spaces, for example, appear with canonical orbifold structures.

   Recently Chen and Ruan \cite{ChenRuan} motivated by  their ideas in quantum cohomology
and by orbifold string theory models discovered a remarkable cohomology theory of
orbifolds that they have coined \emph{orbifold cohomology}. Adem and Ruan
\cite{AdemRuan} went on to define the corresponding \emph{orbifold K-theory} and to
study the resulting Chern isomorphism. One of the remarkable properties of the theory
is that both theories can be twisted by what Ruan has called an \emph{inner local
system} coming from a third group-cohomology class called \emph{discrete torsion}.

   Independently of that, Witten \cite{Witten}, while studying $K$-theory as the natural
recipient of the charge of a D-brane in type IIA superstring
theories was motivated to define a twisting of $K(M)$ for $M$
smooth by a third cohomology class in $H^3(M)$ coming from a
codimension 3-cycle in $M$ and Poincar\'e duality. This twisting
appeared previously in the literature in different forms
\cite{DonovanKaroubi, GrothendieckBrauer, Rosenberg}.

   In this paper we show that if an orbifold
is interpreted as a \emph{stack} then we can define a twisting of
the natural $K$-theory of the stack that generalizes both,
Witten's and Adem-Ruan's twistings. We also show how we can
interpret the theory of \emph{bundle gerbes} over a smooth
manifold and their $K$-theory \cite{Murray0, Murray1, Murray2} in
terms of the theory developed here.

Since the approach to the theory of stacks that we will follow is
not yet published \cite{Fulton}, we try very hard to work in very
concrete terms and so our study includes a very simple definition
of a gerbe over a stack motivated by that of Chaterjee and Hitchin
\cite{Hitchin} on a smooth manifold. This definition is easy to
understand from the point of view of differential geometry, and
of algebraic geometry.

   Using results of Segal \cite{Segal1, Segal2} on
the topology of classifying spaces of categories and of Crainic,
Moerdijk and Pronk on sheaf cohomology over orbifolds
\cite{CrainicMoerdijk, MoerdijkPronk1, MoerdijkPronk, Moerdijk} we
show that the usual theory for the characteristic class of a gerbe
over a smooth manifold \cite{Brylinski} extends to the orbifold
case. Then we explain how Witten's arguments relating the charge
of a $D$-brane generalize.

   A lot of what we will show is valid for foliation groupoids and also for
    a category of \emph{Artin stacks} - roughly speaking
spaces that are like orbifolds except that we allow the stabilizers of the local
actions to be Lie groups. In particular we will explain how the twisting proposed
here can be used to realize the Freed-Hopkins-Teleman twisting used in their
topological interpretation of the Verlinde algebra \cite{Freed}.

   We would like to use this opportunity to thank Alejandro Adem, Marius Crainic, Bill Dwyer,
Paulo Lima-Filho, Bill Fulton, Mike Hopkins, Shengda Hu, Haynes Miller, Ieke Moerdijk, Mainak Poddar,
Joel Robbin and Yong-Bin Ruan for very helpful discussions regarding this work that
grew out of seminars at Michigan by Fulton, and at Wisconsin by Adem and Ruan. We
would also thank letters by Jouko Mickelsson, 
Michael Murray, Eric Sharpe, Zoran Skoda and very specially Ieke Moerdijk and Alan
Weinstein regarding the preliminary version of this paper.

\section{A review of Orbifolds}

In this section we will review the classical construction of the category of
orbifolds. This category of orbifolds is essentially that introduced by Satake
\cite{Satake} under the name of \emph{V-manifolds}, but with a fundamental difference
introduced by Chen and Ruan. They have restricted the morphisms of the category from
orbifold maps to \emph{good maps}, in fact Moerdijk and Pronk have found
this category previously \cite{MoerdijkPronk1} where good maps go by the name of \emph{strict maps}.
This is the correct class of morphisms from the
point of view of stack theory as we will see later.

\subsection{Orbifolds, good maps and orbibundles} \label{sectionorbifolds}

Following Ruan \cite{Ruan, ChenRuan} we will use the following definition for an
\emph{orbifold}.

\begin{definition}
An n-dimensional \emph{uniformizing system} for a connected
topological space $U$ is a triple $(V,G,\pi)$ where
\begin{itemize}
\item $V$ is a connected $n$-dimensional smooth manifold
\item $G$ is a finite group acting on $V$ smoothly ($C^\infty$ automorphisms)
\item $\pi : V \longrightarrow U$ is a continuous map inducing a homeomorphism
$\tilde{\pi}: V/G \to U$
\end{itemize}
\end{definition}

Two uniformizing systems, $(V_1,G,_1\pi_1)$ and $(V_2,G_2,\pi_2)$
are \emph{ isomorphic} if there exists a pair of functions
$(\phi,\lambda)$ such that:

\begin{itemize}
\item $\phi : V_1 \longrightarrow V_2$ a diffeomorphism
\item $\lambda: G_1 \longrightarrow G_2$ an isomorphism
\end{itemize}

with $\phi$ being $\lambda$-equivariant and $\pi_2 \circ \phi = \pi_1$.

Let $i : U' \hookrightarrow U$ a connected open subset of $U$ and $(V',G',\pi')$
a uniformizing system of $U'$.

\begin{definition} \label{induced}
$(V',G',\pi')$ is \emph{induced} from $(V,G,\pi)$ if there exist:

\begin{itemize}
\item a monomorphism $\lambda: G' \to G$ inducing an isomorphism $\lambda: ker \ G' 
\stackrel{\cong}{\to} ker \ G$, where $ker \ G'$ and $ker \ G$ are the subgroups of $G'$ and $G$ respectively
that act trivially on $V'$ and $V$, and
\item a $\lambda$-equivariant open embedding $\phi: V' \to V$
\end{itemize}

with $i \circ \pi' = \pi \circ \phi$. We call $(\phi,\lambda) :
(V',G',\pi') \to (V,G,\pi)$ an \emph{injection}.
\end{definition}

Two injections $(\phi_i,\lambda_i) : (V'_i,G'_i,\pi'_i) \to (V,G,\pi)$, $i=1,2$, are
 isomorphic if there exist:

\begin{itemize}
\item an isomorphism $(\psi,\tau) : (V'_1,G'_1,\pi'_1) \to (V'_2,G'_2,\pi'_2)$ and
\item an automorphism $(\tilde{\psi},\tilde{\tau}) : (V,G,\pi) \to (V,G,\pi)$
\end{itemize}
such that $(\tilde{\psi},\tilde{\tau}) \circ (\phi_1,\lambda_1) = (\phi_2,\lambda_2) \circ (\psi,\tau)$

\begin{rem} \label{induceduniformizingsystem}
Since for a given uniformizing system $(V,G,\pi)$ of $U$, and any connected open set
$U'$ of $U$, $(V,G,\pi)$ induces a unique
isomorphism class of uniformizing systems of $U'$ 
we can define the \emph{germ} of a uniformizing system localized
at a point.
\end{rem}

Let $U$ be a connected and locally connected topological space, $p \in U$ a point, and
$(V_1,G_1,\pi_1)$ and $(V_2,G_2,\pi_2)$ uniformizing systems of the neighborhoods
$U_1$ and $U_2$ of $p$ respectively, then

\begin{definition}
$(V_1,G_1,\pi_1)$ and $(V_2,G_2,\pi_2)$ are { equivalent at $p$} if they induce
uniformizing systems for a neighborhood $U_3 \subset U_1 \cap U_2$ of $p$
\end{definition}

The \emph{germ} of $(V,G,\pi)$ at $p$ is defined as the set of
uniformizing systems of neighborhoods of $p$ which are equivalent
at $p$ with $(V,G,\pi)$.

\begin{definition} \label{orbifoldstructure}
Let $X$ be a Hausdorff, second countable topological space. An
{$n$-dimensional orbifold structure} on $X$ is a set
$\{(V_p,G_p,\pi_p) | p \in X\}$ such that
\begin{itemize}
\item $(V_p,G_p,\pi_p)$ is a uniformizing system of $U_p$, neighborhood of $p$ in $X$
\item for any point $q \in U_p$,  $(V_p,G_p,\pi_p)$ and $(V_q,G_q,\pi_q)$ are equivalent
at $q$
\end{itemize}
\end{definition}

We say that two orbifold structures on $X$, $\{(V_p,G_p,\pi_p)\}_{p \in X}$ and
$\{(V'_p,G'_p,\pi'_p)\}_{p \in X}$, are {equivalent} if for any $q \in X$
$(V_q,G_q,\pi_q)$ and $(V'_q,G'_q,\pi'_q)$ are equivalent at $q$.

\begin{definition}
With a given orbifold structure, $X$ is called an  orbifold.
\end{definition}

Sometimes we will simply denote by $X$ the pair
$(X,\{(V_p,G_p,\pi_p)\}_{p \in X})$. When we want to make the
distinction between the underlying topological space $X$ and the
orbifold $(X,\{(V_p,G_p,\pi_p)\}_{p \in X})$ we will write $\Xx$
for the latter.

For any $p \in X$ let $(V,G,\pi)$ be a uniformizing of a neighborhood around $p$ and 
$\bar{p} \in \pi^{-1}(x)$. Let $G_p$ be the stabilizer of $G$ at $p$. Up to conjugation
the group $G_p$ is independent of the choice of $\bar{p}$ and is called the {\emph{ 
isotropy group}} or {\emph{ local group}} at $p$.

\begin{definition}
An orbifold $X$ is called {\emph{reduced}} if the isotropy groups $G_p$ act effectively
for all $p \in X$.
\end{definition}
In particular this implies that an orbifold is reduced if and only if the groups
$ker \ G$ of definition \ref{induced}  are all trivial.

\begin{example} Let $X=Y/G$ be the orbifold which is the global quotient of the finite
group $G$ acting on a connected space $Y$ via automorphisms. Then
$\{(X,G,\pi)\}$ is trivially an orbifold structure for $X$. We
can also define another equivalent orbifold structure in the
following way: for $p \in X$, let $U_p \subset X$ be a
sufficiently small neighborhood of $p$ such that $$\pi^{-1}(U_p)
= \bigsqcup_\alpha V_p^\alpha$$ the disjoint union of
neighborhoods $V_p^\alpha$, where $G$ acts as permutations on the
connected components of $\pi^{-1}(U_p)$.

Let $V_p$ be one of these connected components, and let $G_p$ be the subgroup of $G$
which fixes this component $V_p$ (we could have taken $U_p$ so that $G_p$ is the isotropy
group of the point $y \in \pi^{-1}(p) \cap V_p$) and take $\pi_p = \pi|_{V_p}$, then $V_p /G_p
\stackrel{\cong}{\to} U_p$ and $(V_p, G_p, \pi_p)$ is a uniformizing system for $U_p$. This
is a direct application of the previous remark.
\end{example}

Now we can define the notion of an orbifold vector bundle or
orbibundle of rank $k$. Given a uniformized topological space $U$
and a topological space $E$ with a surjective continuous map $pr:
E \to U$, a uniformizing system of a rank $k$ vector bundle $E$
over $U$ is given by the following information:
\begin{itemize}
\item A uniformizing system $(V,G,\pi)$ of $U$
\item A uniformizing system $(V \times \real^k, G, \tilde{\pi})$ for $E$ such that the action of $G$ on $V \times \real^k$
is an extension of the action of $G$ on $V$ given by $g(x,v) = (gx, \rho(x,g)v)$ where
$\rho : V \times G \to Aut(\real^k)$ is a smooth map which satisfies
$$\rho(gx,v) \circ \rho(x,g) = \rho(x, h \circ g) , \ \ g,h \in G, x \in V$$
\item The natural projection map $\tilde{pr} :V \times \real^k \to V$ satisfies
$\pi \circ \tilde{pr} = pr \circ \tilde{\pi}$.
\end{itemize}

In the same way the orbifolds were defined, once we have the uniformizing systems of rank $k$ we can define the germ of
orbibundle structures.

\begin{definition} \label{definitionorbibundle}
 The topological space $E$ provided with a given germ of vector bundle structures over the
orbifold structure of $X$, is an orbibundle over $X$.
\end{definition}

Let's consider now orbifolds $X$ and $X'$ and a continuous map $f:X \to X'$. A
\emph{lifting} of $f$ is the following: for any point $p \in X$ there are charts
$(V_p,G_p,\pi_p)$ at $p$ and $(V_{f(p)},G_{f(p)},\pi_{f(p)})$
 at $f(p)$, and a lifting $\tilde{f}_p$ of $f_{\pi_p(V_p)} :\pi_p(V_p) \to \pi_{f(p)}(V_{f(p)})$ such that for any
$q \in \pi_{p}(V_p)$, $\tilde{f}_q$ and $\tilde{f}_p$ define the same germ of liftings of $f$ at $q$.
\begin{definition}
A $C^\infty$ map between orbifolds  $X$ and $X'$ (orbifold-map) is a germ of $C^\infty$ liftings of a continuous map
between $X$ and $X'$.
\end{definition}
We would like to be able to pull-back bundles using maps between
orbifolds, but it turns out that with general orbifold-maps they
cannot be defined. We need to restrict ourselves to a more
specific kind of maps between orbifolds, they were named
\emph{good maps} by Chen and Ruan (see \cite{ChenRuan}). These
good maps will precisely match the definition of a morphism in
the category of groupoids (see Proposition
\ref{morphism=goodmap}).

Let $\tilde{f}: X \to X'$ be a $C^\infty$ orbifold-map
 whose underlying continuous function is
$f$. Suppose there is a compatible cover $\UU$ of $X$ and a collection of open
subsets $\UU'$ of $X'$ defining the same germs, such that there is a $1-1$
correspondence between elements of $\UU$ and $\UU'$, say $U \leftrightarrow U'$,
with $f(U) \subset U'$ and $U_1 \subset U_2$ implies $U_1' \subset U_2'$. Moreover,
there is a collection of local $C^\infty$ liftings of $f$ where
 $\tilde{f}_{UU'}:(V,G,\pi) \to (V', G', \pi,)$ satisfies that for each injection
$(i,\phi): (V_1,G_1,\pi_1) \to (V_2, G_2, \pi_2)$ there is another
injection associated to it $(\nu(i),\nu(\phi)): (V_1',G_1',\pi_1') \to (V_2',
G_2', \pi_2')$ with $\tilde{f}_{U_1U_1'} \circ i = \nu(i) \circ
\tilde{f}_{U_2 U_2'}$; and for any composition of injections $j
\circ i$, $\nu(j \circ i) = \nu(j) \circ \nu(i)$ should hold.

The collection of maps $\{\tilde{f}_{UU'}, \nu \}$ defines a $C^\infty$ lifting of
$f$. If it is  in the same germ as $\tilde{f}$ it is called a \emph{compatible
system} of $\tilde{f}$.

\begin{definition} \label{definitiongoodmap}
A $C^\infty$ map is called \emph{good} if it admits a compatible system.
\end{definition}

\begin{lemma} \cite[Lemma 2.3.2]{Ruan} Let $pr:E \to X$ be and orbifold vector bundle over $X'$. For
any compatible system $\xi= \{\tilde{f}_{UU'}, \nu \}$ of a good $C^\infty$ map $f: X
\to X'$, there is a canonically constructed pull-back bundle of $E$ via $\tilde{f}$
(a bundle $pr :E_\xi \to X$ together with a $C^\infty$ map $\tilde{f}_\xi : E_\xi \to
E$ covering $\tilde{f}$.)
\end{lemma}

\subsection{Orbifold Cohomology} \label{sectionorbifoldcohomology}

Motivated by index theory and by string theory Chen and Ruan have
defined a remarkable cohomology theory for orbifolds. One must
point out that while as a group it had appeared before in the
literature in several forms, its product is completely new and
has very beautiful properties.

For $X$ an orbifold, and $p$ a point in $U_p \subset X$ with
$(V_p, G_p, \pi_p), \pi(V_p) = U_p$
 a local chart around it, the multi-sector
$\widetilde{\Sigma_k X}$ is defined as the set of pairs $(p,
({\bf g)})$, where $({\bf g})$ stands for the conjugacy class of
${\bf g} = (g_1, \dots , g_k)$ in $G_p$.

For the point $(p,({\bf g})) \in \widetilde{\Sigma_k X} $  the multisector can be seen locally as
$$ V_p^{\bf g}/C({\bf g})$$
where $V_p^{\bf g} = V_p^{g_1} \cap V_p^{g_2} \cap \cdots \cap
V_p^{g_k}$ and $ C({\bf g})) = C(g_1) \cap C(g_2) \cap \cdots
\cap C(g_k)$. $V_p^g$ stands for the fixed-point set of
 $g \in G_p$ in $V_p$, and $C(g)$ for the centralizer of $g$ in $G_p$.

Its connected components are described  in the following way.
 For $q \in U_p$, up to conjugation,
there is a injective homomorphism $G_q \to G_p$, so for ${\bf g}$
in $G_q$ the conjugacy class $({\bf g})_{G_p}$ is well defined.
In this way we can define an equivalence relation $({\bf
g})_{G_q} \cong ({\bf g})_{G_p}$ and we call $T_k$ the set of
such equivalence classes. We will abuse of notation and will
write $({\bf g})$ to denote the equivalence class at which $({\bf
g})_{G_p}$ belongs to. Let $T_k^0 \subset T_k$ be the set of
equivalence classes $({\bf g})$ such that $g_1 g_2 \dots g_k =1$.

$\widetilde{\Sigma_k X}$ is decomposed as $\bigsqcup_{({\bf g})
\in T_k} X_{({\bf g})}$ where
$$X_{({\bf g})} = \left\{ (p,({\bf g'})_{G_p}) | {\bf g'} \in G_p \ \  \& \ \ ({\bf g'})_{G_p} \in ({\bf g}) \right\}$$
$X_{(g)}$ for $g \neq 1$ is called a twisted sector and $X_{(1)}$
the non-twisted one.

\begin{example}
Let's consider the global quotient $X=Y/G$, $G$ a finite group. Then $X_{(g)} \cong Y^g/C(g)$ where $Y^g$
is the fixed-point set of $g \in G$ and $C(g)$ is its centralizer, hence
$$\widetilde{\Sigma X} \cong \bigsqcup_{\{(g) | g \in G \}} Y^g /C(g) $$
\end{example}

Let's consider the natural maps between multi-sectors; the
evaluation maps $e_{i_1, \dots, i_l} : \widetilde{\Sigma_k X} \to
\widetilde{\Sigma_l X}$ defined by $e_{i_1, \dots,
i_l}(x,(g_1,\dots, g_k)) \mapsto (x,(g_{i_1}, \dots, g_{i_l}))$
and the involutions $I:\widetilde{\Sigma_k X} \to
\widetilde{\Sigma_k X}$ defined by $I(x,({\bf g})) \mapsto
(x,({\bf g^{-1}}))$ where ${\bf
g^{_1}}=(g_1^{-1},\dots,g_k^{-1})$.

An important concept in the theory is that of an \emph{inner local system} as defined
by Y. Ruan \cite{Ruan}. We will show below that inner local systems are precisely
modeled by \emph{gerbes} over the orbifold.

\begin{definition}
 \label{innerlocalsystem}
Let $X$ be an orbifold. An inner local system
$\LL=\{L_{(g)}\}_{(g) \in T_1}$  is an assignment of a flat
complex line orbibundle $L_{(g)} \to X_{(g)}$ to each twisted
sector $X_{(g)}$ satisfying the compatibility conditions:
\begin{enumerate}
\item $L_{(1)}= 1$ is trivial.
\item $I^*L_{(g^{-1})} = L_{(g)}$
\item Over each $X_{({\bf g})}$ with $({\bf g}) \in T_3^0$ ($g_1g_2g_3=1$),
$$e_1^*L_{(g_1)} \otimes e_2^*L_{(g_2)} \otimes e_3^*L_{(g_3)} = 1 $$
\end{enumerate}
\end{definition}
One way to introduce inner local systems is by discrete torsion.
Let $Y$ be the universal orbifold cover of the orbifold $Z$, and
let $\pi_1^{orb}(Z)$  be the group of deck transformations (see
\cite{Thurston}).

For $X=Z/G$, $Y$ is an orbifold universal cover of $X$ and we
have the following short exact sequence:
$$ 1 \To \pi_1(Z) \To \pi_1^{orb}(X) \To G \To 1$$

We call an element in $ H^2(\pi_1^{orb}(X), \U{1})$ a discrete
torsion of $X$. Using the previous short exact sequence
$H^2(G,\U{1}) \to H^2(\pi_1^{orb}(X), \U{1})$, therefore elements
$\alpha \in H^2(G,\U{1})$ induce
 discrete torsions.

We can see $\alpha : G \times G \to \U{1}$ as a two-cocycle
satisfying $\alpha_{1,g} = \alpha_{g,1} =1 $ and $\alpha_{g,hk}
\alpha_{h,k} = \alpha_{g,h} \alpha_{gh,k}$ for any $g,h,k \in G$.
We can define its phase as $\gamma(\alpha)_{g,h} := \alpha_{g,h}
\alpha_{h,g}^{-1}$ which induces a representation of $C(g)$
$$L^\alpha_g := \gamma(\alpha)_{g,\_}: C(g) \to \U{1}$$

\begin{example}
In the case that $Y \to X$ is the orbifold universal cover and
$G$ is  the orbifold fundamental group such that $X=Y/G$, we can
construct a complex line bundle $L_g = Y^g \times_{L^\alpha_g}
\complex$ over $X_{(g)}$. We get that $L_{tgt^{-1}}$ is naturally
isomorphic to $L_g$ so we can denote the latter one by $L_{(g)}$,
and restricting to $X_{(g_1,\dots,g_k)}$, $L_{(g_1, \dots,
g_k)}=L_{(g_1)} \cdots L_{(g_k)}$; then $\LL= \{L_{(g)} \}_{(g)
\in T_1}$ is an inner local system for $X$
\end{example}

To define the orbifold cohomology group we need to add a shifting to  the cohomology
of the twisted sectors, and for that we are going to assume that the orbifold $X$ is
almost complex with complex structure $J$; recall that $J$ will be a smooth section
of $End(TX)$ such that $J^2 =-Id$.

For $p \in X$ the almost complex structure gives rise to an effective representation
$\rho_p : G_p \to \GLC{n}$ ($n = dim_{\complex} X$) that could be diagonalized as
$$diag \left(e^{2 \pi \frac{m_{1,g}}{m_g}}, \dots , e^{2 \pi \frac{m_{n,g}}{m_g}} \right)$$
where $m_g$ is the order of $g$ in $G_p$ and $0 \leq m_{j,g} < m_g$. We define a
function $\iota : \widetilde{ \Sigma X} \to \rational$ by
$$\iota(p,(g)_{G_p}) = \sum_{j=1}^n \frac{m_{j,g}}{m_g}$$
It is easy to see that it is locally constant, hence we call it
$\iota_{(g)}$; it is an integer if and only if $\rho_p(g) \in
\SLC{n}$ and
$$\iota_{(g)} + \iota_{(g^{-1})} = rank(\rho_p(g) - I)$$
which is the complex codimension $dim_{\complex} X - dim_{\complex} X_{(g)}$.
$\iota_{(g)}$ is called the degree shifting number.

\begin{definition}
Let $\LL$ be an inner local system, the orbifold cohomology groups are defined as
$$H_{orb}^d(X; \LL) = \bigoplus_{(g) \in T_1} H^{d- 2 \iota_{(g)}} (X_{(g)}; \LL_{(g)})$$
If $\LL = \LL_\alpha$ for some discrete torsion $\alpha$ we define
$$H^*_{orb,\alpha} (X, \complex) = H^*_{orb} (X, \LL_\alpha)$$
\end{definition}

\begin{example}
For the global quotient $X=Y/G$ and $\alpha \in H^2(G, \U{1})$, $L_g^\alpha$ induces
a twisted action of $C(g)$ on the cohomology of the fixed point set $H^*(Y^g,
\complex)$ by $\beta \mapsto L_g^\alpha(h) h^* \beta$ for $h \in C(g)$. Let $H^*(Y^g,
\complex)^{C^\alpha (g)}$ be the invariant subspace under this twisted action. Then
$$H^d_{orb, \alpha} (X ; \complex) = \bigoplus_{(g) \in T_1} H^{d - 2 \iota_{(g)}} (Y^g ; \complex)^{C^\alpha(g)}$$
\end{example}

\subsection{Orbifold $K$-theory} \label{subsectionorbifoldKtheory}

In this section we will briefly describe a construction by Adem and Ruan of the
so-called twisted orbifold $K$-theory. The following construction will generate a
twisting of the orbifold $K$-theory by certain class in a group cohomology group. We
will recover this twisting later, as a particular case of a twisting of $K$-theory on
a groupoid by an arbitrary gerbe.

The following constructions are based on projective representations. A function $\rho
: G \to GL(V)$, for $V$ a finite dimensional complex vector space, is a projective
representation of $G$ if there exists a function $\alpha :G \times G \to \complex^*$
such that $\rho(x) \rho(y) = \alpha(x,y) \rho(xy)$. Such $\alpha$ defines a
two-cocycle on $G$, and $\rho$ is said to be $\alpha$-representation on the space
$V$. We can take sum of any two $\alpha$-representations, hence we can define the
Grothendieck group, $R_\alpha(G)$
 associated to the monoid of linear isomorphism classes of such $\alpha$-representations.

Let's assume that $\Gamma$ is a semi-direct product of a compact Lie group $H$ and
a discrete group $G$. Let $\alpha \in H^2(G,\U{1})$ so we have a group extension
$$1 \to \U{1} \to \tilde{G} \to G \to 1$$
and $\tilde{\Gamma}$ is the semi-direct product of $H$ and $\tilde{G}$.

Suppose that $\Gamma$ acts on a smooth manifold $X$ such that $X/\Gamma$ is compact
and the action has only finite isotropy, then $Y=X/\Gamma$ is an orbifold.

\begin{definition}
An $\alpha$ twisted $\Gamma$-vector bundle on $X$ is a complex
vector bundle $E \to X$ such that $\U{1}$ acts on the fibers
through complex multiplication extending the action of $\Gamma$
in $X$ by an action of $\tilde{\Gamma}$ in $E$.

We define ${}^\alpha K_\Gamma (X)$ the $\alpha$-twisted $\Gamma$-equivariant
$K$-theory of $X$ as the Grothendieck group of isomorphism classes of $\alpha$ twisted
$\Gamma$-bundles over $X$.
\end{definition}

For an $\alpha$-twisted bundle $E \to X$ and a $\beta$-twisted
bundle $F \to X$ consider the tensor product bundle $E \otimes F
\to X$, it becomes an $\alpha +\beta$-twisted bundle. So we have
a product
$${}^\alpha K_\Gamma (X) \otimes {}^\beta K_\Gamma(X) \to {}^{\alpha + \beta} K_\Gamma(X)$$
And so we call the \emph{total twisted equivariant K-theory of a $\Gamma$ space} as:
$$TK_\Gamma(X) = \bigoplus_{\alpha \in H^2(G,\U{1})} {}^\alpha K_\Gamma (X)$$

When $\Gamma$ is a finite group, there is the following decomposition theorem,

\begin{theorem}\cite[Th. 4.2.6]{Ruan}\label{ADEMRUANTH}
Let  $\Gamma$ be a finite group that acts on $X$,
 then for any $\alpha \in H^2(G,\U{1})$
$${}^\alpha K^*_\Gamma(X) \otimes \complex \cong H^*_{orb,\alpha}(X/\Gamma; \complex)$$
\end{theorem}

The decomposition is as follows:
$${}^\alpha K^*_\Gamma(X) \otimes \complex \cong \bigoplus_{(g)} (K(X^g) \otimes L^\alpha_g)^{C(g)}
\cong \bigoplus_{(g)} H^*(X^g ; \complex)^{C^\alpha (g)} \cong
H^*_{orb,\alpha}(X/\Gamma; \complex)$$

\begin{definition}\label{alfatwistedk} In the case that $Y \to X$ is the orbifold universal cover and $\alpha \in
H^2(\pi_1^{orb}(X), \U{1})$, the $\alpha$ twisted orbifold $K$-theory, ${}^\alpha
K_{orb}(X)$, is the Grothendieck group of isomorphism classes of $\alpha$-twisted
$\pi_1^{orb}(X)$-orbifold bundles over $Y$ and the total orbifold $K$-theory is:
$$TK_{orb}(X) = \bigoplus_{\alpha \in H^2(\pi_1^{orb}(X); \U{1})} {}^\alpha K_{orb}(X)$$
\end{definition}

\subsection{Twisted $K$-theory on smooth
manifolds.}\label{smoothtwist}

In \cite{Witten} Witten shows that the $D$-brane charge for Type
IIB superstring theories (in the case of 9-branes) should lie on
a twisted $K$-theory group that he denotes as $K_{[H]}(X)$ where
a 3-form $H\in\Omega^3(X;\real)$ models the Neveu-Schwarz
$B$-field and $[H]\in H^3(X;\integer)$ is an integer cohomology
class. The manifold $X$ is supposed smooth and it is where the
$D$-branes can be wrapped. The class $[H]$ is not torsion, but in
any case when $[H]$ is a torsion class Witten gives a very
elementary definition of $K_{[H]}(X)$. This will also be a
particular class of the twisting of K-theory on a stack by a
gerbe defined below.

The construction of $K_{[H]}(X)$ is as follows. Consider the long exact sequence in
\emph{simplicial} cohomology with \emph{constant} coefficients
\begin{equation}
   \cdots\to H^2(X;\real)\toparrow{i} H^2(X;\U{1})\to
   H^3(X;\integer)\to H^3(X;\real)\to\cdots
\end{equation}
induced by the exponential sequence
$0\to\integer\toparrow{i}\real\toparrow{\exp}\U{1}\to1$. Since $[H]$ is torsion, it
can be lifted to a class $H^*\in H^2(X;\U{1})$, and if $n$ is its order then for a
fine covering $\UU=\{U_i\}_i$ of $X$ the class $H^*$ will be represented by a
\v{C}ech cocycle $h_{ijk} \in \check{C}^3(X)(\rational(\zeta_n))$ valued on $n$-th
roots of unity.

Now we can consider a vector bundle as a collection of functions
$g_{ij}\colon U_{ij} \to \GLC{m}$ such that
$g_{ij}g_{jk}g_{ki}=id_{\GLC{m}}$.

\begin{definition}\label{Ktwistedsmooth}
We say that a collection of functions $g_{ij}\colon U_{ij} \to
\GLC{m}$ is an \emph{$[H]$-twisted vector bundle} $E$ if
$g_{ij}g_{jk}g_{ki}=h_{ijk} \cdot id_{\GLC{m}}$. The Grothendieck group of such
twisted bundles is $K_{[H]}(X)$.
\end{definition}

This definition does not depend on the choice of cover, for it can be written in terms
of a Grothendieck group of modules over the algebra of sections ${\mathrm{END}}(E)$
of the endomorphism bundle $E\otimes E^*$, that in particular \emph{is an ordinary
vector bundle} \cite{DonovanKaroubi}.

In the case in which the class $\alpha=[H]$ is not a torsion class one can still
define a twisting and interpret it in terms of Fredholm operators on a Hilbert space.
The following description is due to Atiyah and Segal \cite{Atiyah}. Let $\HH$ be a
fixed Hilbert space. We let $\BB(\HH)$ be the Banach algebra of bounded operators on
$\HH$ and $\FF(\HH)\subset\BB(\HH)$ be the space of Fredholm operators on $\HH$,
namely, those operators in $\BB(\HH)$ that are invertible in $\BB(\HH)/\KK(\HH)$
where $\KK(\HH)$ is the ideal in $\BB(\HH)$ consisting of compact operators.

 Then we have the following classical theorem of Atiyah,
\begin{equation*}
   K(X)=[X,\FF],
\end{equation*}
where the right hand side means all the homotopy classes of maps
$X\to\FF$. In particular $\FF\simeq BU$.

   For a cohomology class $\alpha \in H^3(X, \integer)$
Atiyah and Segal construct a bundle $\FF_\alpha$ over $X$ with
fiber $\FF(\HH)$, and then define the twisted $K_\alpha$-theory as
\begin{equation}\label{AtiyahTwist}
 K_\alpha (X) =[\Gamma(\FF_\alpha)]
\end{equation} namely the homotopy classes of sections of the bundle
$\FF_\alpha$.

To construct $\FF_\alpha$ notice that it is enough to construct a bundle $\BB_\alpha$
over $X$ with fiber $\BB(\HH)$ for  a given class $\alpha$. Observe now the property
that a bounded linear map $\HH\to\HH$ is Fredholm is completely determined by the map
$\Proj(\HH) \to \Proj(\HH)$ that it induces. So it will be enough to construct a
infinite dimensional projective bundle $\Proj_\alpha$ with fiber $\Proj(\HH)$. This
can be done using Kuiper's theorem that states that the group $U(\HH)$ of unitary
operators in $\HH$ is contractible and therefore one has
$\Proj(\complex^\infty)=K(\integer,2)=BU(1)=U(\HH)/U(1)=\Proj U(\HH)$ and
$K(\integer,3)\simeq B\Proj U(\HH)$. The class $\alpha \in
H^3(X,\integer)=[X,K(\integer,3)]=[X,B\Proj U(\HH)]$ gives the desired projective
bundle, at it is called the \emph{Dixmier-Douady} class of the projective bundle. It
is worthwhile to mention that J. Rosenberg has previously defined $K_\alpha (X)$ in
\cite{Rosenberg}. His definition is clearly equivalent to the one explained above.

\section{Gerbes over smooth manifolds} \label{sectiongerbesmanifolds}

\subsection{Gerbes}\label{gerbessmooth}
As a way of motivation for what follows, later we will summarize the facts about
gerbes over smooth manifolds, we recommend to see \cite{Brylinski, Hitchin} for a
more detailed description of the subject. Just as a line bundle can be given by
transition functions, a gerbe can be given by transition data, namely line bundles.
But the ``total space" of a gerbe is a stack, as explained in the appendix. The same
gerbe can be given as transition data in several ways.

Let's  suppose $X$ is a smooth manifold and $\{U_\alpha
\}_\alpha$ an open cover. Let's consider the functions
$$g_{\alpha \beta \gamma} : U_\alpha \cap U_\beta \cap U_\gamma \To \U{1}$$
defined on the threefold intersections satisfying
$$g_{\alpha \beta \gamma} = g_{\alpha \gamma \beta}^{-1} = g_{\beta \alpha \gamma}^{-1} = g_{\gamma \beta \alpha}^{-1}$$
and the cocycle condition
$$(\delta g)_{\alpha \beta \gamma \eta} = g_{\beta \gamma \eta} g_{\alpha \gamma \eta}^{-1}
g_{\alpha \beta \eta} g_{\alpha \beta \gamma}^{-1} =1 $$ on the four-fold
intersection $U_\alpha \cap U_\beta \cap U_\gamma \cap U_\eta$. All these data define
a \emph{gerbe}. We could think of $g$ as a $\check{C}ech$ cocycle of
$H^2(X,C^\infty(\U{1}))$ and therefore we can tensor them using the product of
cocycles. It also defines a class in $H^3(X; \integer)$; taking the long exact
sequence of cohomology
$$\cdots \to H^i(X,C^\infty(\real)) \to H^i(X,C^\infty(\U{1})) \to H^{i+1}(X, \integer) \to \cdots$$
 given from the exact sequence of sheaves
$$ 0 \to \integer \to C^\infty(\real) \to C^\infty(\U{1}) \to 1$$
and using that $C^\infty(\real)$ is a fine sheaf,  we get $H^2(X, C^\infty(\U{1}))
\cong H^3(X,\integer)$. We might say that a gerbe is determined topologically by its
characteristic class.

A trivialization of a gerbe is defined by functions
$$f_{\alpha \beta} =f_{\beta \alpha} : U_\alpha \cap U_\beta \to \U{1}$$
on the twofold intersections such that
 $$g_{\alpha \beta \gamma} = f_{\alpha \beta} f_{\beta \gamma} f_{\gamma \alpha}$$
In other words $g$ is represented as a coboundary $\delta f = g$.

The difference  of two trivializations $f_{\alpha \beta}$ and
$f'_{\alpha \beta}$ given by $h_{\alpha \beta}$ becomes a line
bundle ($h_{\alpha \beta} h_{\beta \gamma} h_{\gamma \alpha} =1$).

Over  a particular open subset $U_0$ we can define a
trivialization, for $\beta, \gamma \neq 0$ we take $f_{\beta
\gamma} := g_{0 \beta \gamma}$ and because of  the cocycle
condition we have $g_{\alpha \beta \gamma} = f_{\alpha \beta}
f_{\beta \gamma} f_{\gamma \alpha}$. Adding $f_{0 \beta} = 1$ we
get a trivialization localized at $U_0$ and we could do the same
over each $U_\alpha$. Then on the intersections $U_\alpha \cap
U_\beta$ we get two trivializations that differ by a line bundle
$\LL_{\alpha \beta}$. Thus a gerbe can also be seen as the
following data:

\begin{itemize}
\item A line bundle $\LL_{\alpha \beta}$ over each $U_\alpha \cap U_\beta$
\item $\LL_{\alpha \beta} \cong \LL_{\beta \alpha}^{-1}$
\item A trivialization $\theta_{\alpha \beta \gamma}$ of $\LL_{\alpha \beta} \LL_{\beta \gamma} \LL_{\gamma \alpha}\cong 1$ where $\theta_{\alpha \beta \gamma}: U_{\alpha \beta \gamma} \to \U{1}$ is a 2-cocycle.
\end{itemize}

\begin{example} \cite[Ex. 1.3]{Hitchin} \label{Poincare}
Let $M^{n-3} \subset X^n$  be an oriented codimension 3
submanifold of a compact oriented one $X$. Take coordinate
neighborhoods $U_\alpha$ of $X$ along $M$, we could think of them
as $U_\alpha \cong (U_\alpha \cap M) \times \real^3$, and let
$U_0 = X \backslash N(M)$, where $N(M)$ is the closure of a small
neighborhood of $M$, diffeomorphic to the disc bundle in the
normal bundle. We have $$U_0 \cap U_\alpha \cong U_\alpha \cap M
\times \{x \in \real^3 : ||x|| > \epsilon \}$$ and let define the
bundle $\LL_{\alpha 0}$ as the pullback by $x \mapsto x/||x||$ of
the line bundle of degree 1 over $S^2$.

 The line bundles $\LL_{\alpha \beta}= \LL_{\alpha 0} \LL_{0 \beta}^{-1}$ are defined on
$(U_\alpha \cap U_\beta \cap M) \times \{ x \in \real^3 : ||x|| >
\epsilon \}$ and by construction $c_1(\LL_{\alpha \beta}) =0$
over $S^2$, then  they can be extended to trivial ones on the
whole $U_\alpha \cap U_\beta$. This information provides us with
a gerbe and the characteristic class of it in $H^3(X, \integer)$
is precisely the Poincar\'e dual to the homology class of the
submanifold $M$. This is the gerbe that we will use to recover
Witten's twisting of $K$-theory.
\end{example}

\subsection{Connections over gerbes}

We can also do differential geometry over gerbes \cite{Hitchin} and we will describe
what is a  connection over a gerbe.

For $\{U_\alpha \}$ a cover such that all the intersection are contractible (a {\it Leray} cover),
 a  connection will consist of 1-forms over the double intersections $A_{\alpha \beta}$, such that
$$iA_{\alpha \beta} + i A_{\beta \gamma} + i A_{\gamma \alpha} = g_{\alpha \beta \gamma}^{-1} d g_{\alpha \beta \gamma}$$
where $g_{\alpha \beta \gamma}: U_\alpha \cap U_\beta \cap U_\gamma \to \U{1}$ is the cocycle defined by the gerbe.

Because $d (g_{\alpha \beta \gamma}^{-1} d g_{\alpha \beta \gamma}) =0$ there are 2-forms $F_\alpha$ defined
over $U_\alpha$ such that $F_\alpha - F_\beta = dA_{\alpha \beta}$; as $dF_\alpha = dF_\beta$
 then we define a global 3-form $G$ such that $G|_{U_\alpha} = dF_\alpha$. This 3-form $G$ is called the curvature of
the gerbe connection.

As  the $A_{\alpha \beta}$ are 1-forms over the double
intersections, we could reinterpret them as connection forms over
the line bundles. So, using the line bundle definition of gerbe,
a connection in that formalism is:

\begin{itemize}
\item A connection $\Delta_{\alpha \beta}$ on $\LL_{\alpha \beta}$ such that
\item $\Delta_{\alpha \beta \gamma} \theta_{\alpha \beta \gamma} = 0$ where $\Delta_{\alpha \beta \gamma}$
is the connection over $\LL_{\alpha \beta} \LL_{\beta \gamma} \LL_{\gamma \alpha}$ induced by the
 $\Delta_{\alpha \beta}$
\item A 2-form $F_\alpha \in \Omega^2(U_\alpha)$ such that on $U_\alpha \cap U_\beta$, $F_\beta - F_\alpha$ equals
the curvature of $\Delta_{\alpha \beta}$
\end{itemize}

When the curvature $G$ vanishes we say that the connection on the gerbe is \emph{flat}.

\section{Groupoids}\label{Groupoids}

The underlying idea of everything we do here is that an orbifold is best understood
as a stack. A stack $\XX$ is a ``space" in which we can't talk of a point in $\XX$
but rather only of functions $S\to\XX$ where $S$ is any space, much in the same
manner in which it makes no sense to talk of the value of the Dirac delta $\delta(x)$
at a particular point, but it makes perfect sense to write $\int_\real \delta(x)f(x)
dx$. To be fair there are points in a stack, but they carry automorphism groups in a
completely analogous way to an orbifold. We refer the reader to the Appendix for more
on this. In any case, just as a smooth manifold is completely determined by an open
cover and the corresponding gluing maps, in the same manner a stack will be
completely determined by a groupoid representing it. Of course there may be more than
one such groupoid, so we use the notion of Morita equivalence to deal with this issue.

A groupoid can be thought of as a generalization of a group, a manifold and an
equivalence relation. First an equivalence relation. A groupoid has a set of relations
$\RR$ that we will think of as arrows. These arrows relate elements is a set $\UU$.
Given an arrow $\toparrow{r} \in \RR$ it has a source $x=s(\toparrow{r})\in\UU$ and a
target $y=t(\toparrow{r})\in\UU$. Then we say that $x\toparrow{r} y$, namely $x$ is
related to $y$. We want to have an equivalence relation, for example we want
transitivity and then we will need a way to compose arrows
$x\toparrow{r}y\toparrow{s}z$. We also require $\RR$ and $\UU$ to be more than mere
sets. Sometimes we want them to be locally Hausdorff, paracompact, locally compact
topological spaces, sometimes schemes.

Consider an example. Let $X=S^2$ be the smooth 2-dimensional sphere. Let $p,q$ be the
north and the south poles of $S^2$, and define $U_1=S^2-\{p\}$ and $U_2=S^2-\{q\}$.
Let $U_{12}=U_1\cap U_2$ and $U_{21}=U_2 \cap U_1$ be \emph{two disjoint} annuli.
Similarly take two disjoint disks $U_{11}=U_1\cap U_1$ and $U_{22}=U_2 \cap U_2$.
Consider a category where the objects are $\UU=U_1\uplus U_2$ where $\uplus$ means
disjoint union. The set of arrows will be $\RR=U_{11} \uplus U_{12} \uplus U_{21}
\uplus U_{22}$. For example the point $x\in U_{12} \subset \RR$ is thought  of as an
arrow from $x\in U_1 \subset \UU$ to $x\in U_2 \subset \UU$, namely $x\toparrow{x}
x$. This is a groupoid associated to the sphere. In this example we can write the
disjoint union of all possible triple intersections as $\RR \timests \RR$

\subsection{Definitions}
A \emph{groupoid} is a pair of objects in a category $\RR, \UU$ and morphisms
$$s,t : \RR \rightrightarrows \UU$$ called respectively source and target ,
provided with an identity $$e : \UU \To \RR $$ a multiplication
$$m:\RR \timests \RR \To \RR$$ and an inverse
$$i : \RR \To \RR$$ satisfying the following properties:

\begin{enumerate}
\item The identity inverts both $s$ and $t$:
        $$
        \xymatrix{
        \UU \ar[r]^e \ar[rd]_{id_\UU} & \RR \ar[d]^s\\
          & \UU }
        \ \ \ \ \ \ \ \ \ \ \  \ \ \ \ \
        \xymatrix{
        \UU \ar[r]^e \ar[rd]_{id_\UU} & \RR \ar[d]^t\\
          & \UU }
        $$
\item Multiplication is compatible with both $s$ and $t$:

        $$
        \xymatrix{
        \RR \timests \RR \ar[r]^m \ar[d]_{\pi_1} & \RR \ar[d]^s\\
         \RR \ar[r]^s & \UU }
        \ \ \ \ \ \ \ \ \ \ \  \ \ \ \ \
        \xymatrix{
        \RR \timests \RR \ar[r]^m \ar[d]_{\pi_2} & \RR \ar[d]^t\\
         \RR \ar[r]^t & \UU }
        $$

\item Associativity:

        $$
        \xymatrix{
        \RR \timests \RR \timests \RR \ar[rr]^{id_\RR \times m} \ar[d]_{m \times id_\RR} & & \RR \timests \RR \ar[d]^m\\
         \RR \timests \RR \ar[rr]^m & & \UU }
        $$

\item Unit condition:
        $$
        \xymatrix{
        \RR \ar[r]^(.35){(es, id_\RR)}  \ar[rd]_{id_\RR} & \RR \timests \RR \ar[d]^m\\
          & \RR }
        \ \ \ \ \ \ \ \ \ \ \  \ \ \ \ \
        \xymatrix{
        \RR \ar[r]^(.35){(id_\RR,et)} \ar[rd]_{id_\RR} & \RR \timests \RR \ar[d]^m\\
          & \RR }
        $$

\item Inverse:
$$
\begin{array}{ccc}
i \circ i & = & id_\RR \\
s \circ i & = & t \\
t \circ i & = & s
\end{array}
$$
with

        $$
        \xymatrix{
        \RR  \ar[r]^(.35){(id_\RR, i)} \ar[d]_{s} & \RR \timests \RR \ar[d]^m\\
         \UU \ar[r]^e & \RR }
        \ \ \ \ \ \ \ \ \ \ \  \ \ \ \ \
        \xymatrix{
        \RR  \ar[r]^(.35){(i,id_\RR)} \ar[d]_{t} & \RR \timests \RR \ar[d]^m\\
         \UU \ar[r]^e & \RR }
        $$

\end{enumerate}

We denote the groupoid by  $\RR \twoarrows \UU : = (\RR, \UU,
s,t,e,m,i)$, and the groupoid is called \emph{\'etale } if the
base category is that of locally Hausdorff, paracompact, locally
compact topological spaces and the maps $s,t: \RR \to U$ are
local homeomorphisms (diffeomorphisms). We will say that a
groupoid is proper is $s \times t \colon \RR \to \UU \times \UU$
is a proper (separated) map. We can of course work in the
category of schemes or of differentiable manifolds as well. From
now on we will assume that our groupoids are differentiable,
\'etale and proper; if its also effective, then this groupoid can
be seen as obtained from an orbifold (see \cite[Thm.
4.1]{MoerdijkPronk}).

\medskip
\begin{example} \label{examplemanifold}
For $M$ a manifold and $\{U_\alpha\}$ and open cover, let
$$\UU = \bigsqcup_\alpha U_\alpha \ \ \  \ \RR= \bigsqcup_{(\alpha, \beta)} U_\alpha \cap U_\beta \ \ (\alpha,\beta) \neq (\beta,\alpha)$$
$$s|_{U_{\alpha \beta}}: U_{\alpha \beta} \to U_\alpha, \ t|_{U_{\alpha \beta}}: U_{\alpha \beta} \to U_\beta \ \
e|_{U_{\alpha}}: U_\alpha \to U_\alpha$$
 $$i|_{U_{\alpha \beta}}: U_{\alpha \beta} \to U_{\beta \alpha} \ \& \
m|_{U_{\alpha \beta \gamma}}: U_{\alpha \beta \gamma} \to U_{\alpha \gamma}$$ the
natural maps. Note that in this example $\RR \timests \RR$ coincides with the subset
of $\RR \timests \RR$ of pairs $(u,v)$ so that $t(u)=s(v)$, namely the disjoint union
of all possible triple intersections $U_{\alpha\beta\gamma}$ of open sets in the open
cover $\{U_\alpha\}$. We will denote this groupoid $\RR \rightrightarrows \UU$ by
$\MM(M,U_\alpha)$.
\end{example}

\begin{example}\label{quotient} Let $G$ be a group and $U$ a set provided with a left $G$ action
$$G \times U \To U$$
$$ (g,u) \mapsto gu$$
we put $\UU=U$ and $\RR = G \times U$ with $s(g,u)=u$ and $t(g,u)=gu$. The domain of
$m$ is the same as $G \times G \times U$ where $m(g,h,u)=(gh,u)$,
$i(g,u)=(g^{-1},gu)$ and $e(u)=(id_G,u)$.

We will write $G \times U \twoarrows U$ (or sometimes $\Xx=[U/G]$,) to denote this
groupoid.
\end{example}

\begin{definition}\label{defmor}
A morphism of groupoids $(\Psi, \psi): (\RR' \twoarrows \UU') \To (\RR \twoarrows
\UU)$ are the following commutative diagrams:

        $$\xymatrix{
        \RR' \ar[r]^\Psi \ar@<-.5ex>[d]_{s'} \ar@<+.5ex>[d]^{t'} & \RR \ar@<-.5ex>[d]_{s} \ar@<+.5ex>[d]^t\\
        \UU' \ar[r]^\psi  & \UU}
        \ \ \ \ \ \ \ \ \ \
        \xymatrix{
        \RR' \ar[r]^\Psi  & \RR \\
        \UU' \ar[u]^{e'} \ar[r]^\psi & \UU \ar[u]_e}
        $$
        $$\xymatrix{
        \RR'\: {}_{t'}  \! \times_{s'} \RR' \ar[r]^\Psi \ar[d]_{m'} & \RR \timests \RR \ar[d]^{m} \\
        \RR' \ar[r]^\Psi  & \RR}
        \ \ \ \ \ \ \ \ \ \
        \xymatrix{
        \RR' \ar[r]^\Psi \ar[d]_{i'} & \RR \ar[d]^i\\
        \RR' \ar[r]^\psi  & \RR}
        $$
\end{definition}
Now we need to say when two groupoids are ``equivalent''
\begin{definition}
A morphism of \'{e}tale groupoids $(\Psi,\psi)$ is called an \emph{\'{e}tale Morita morphism} whenever:
\begin{itemize}
\item The map $s \circ \pi_2 : \UU' \fiberprod{\psi}{t} \RR \to \UU$ is an \'etale surjection,
\item The following square is a fibered product
        $$\xymatrix{
        \RR' \ar[r]^\Psi \ar[d]_{(s',t')} & \RR \ar[d]^{(s,t)}\\
        \UU'\times \UU' \ar[r]^{\psi \times \psi}  & \UU \times \UU}
        $$
\end{itemize}
where only the second condition is the required for a morphism of
general groupoids to be Morita. When working on \'{e}tale
groupoids, the Morita morphisms are understood to be \'{e}tale.
\end{definition}

Two groupoids $\RR_1 \twoarrows \UU_1$, $\RR_2 \twoarrows \UU_2$ are called
\emph{Morita equivalent} if there are Morita morphisms $(\Psi_i, \psi_i): \RR'
\twoarrows \UU' \To \RR_i \twoarrows \UU_i$ for $i = 1,2$. This is an equivalence
relation and in general we will consider the category of \'etale groupoids obtained by
formally inverting the Morita equivalences (see \cite{Moerdijk} for details).

It is not hard to define principal bundles over groupoids where the fibers are
groupoids (cf. \cite{CrainicMoerdijk}), but here we will restrict ourselves to the
construction of principal $G$ bundles over groupoids, where $G$ is a Lie group (or an
algebraic group.) This will facilitate the construction of the desired twistings in
$K$-theory.

We give ourselves a groupoid $s,t : \RR \twoarrows \UU$.
\begin{definition} \label{definitionprincipalbundle}
A principal $G$-bundle over the groupoid $\RR \twoarrows \UU$ is the groupoid
 $$\tilde{s}, \tilde{t} : \RR \times G \twoarrows \UU \times G$$ given by the following structure:
$$\tilde{s}(r,h) := (s(r),h)), \ \ \ \ \\ \ \ \tilde{t}(r,h) :=(t(r),\rho(r)h), \ \ \ \ \ \  \ \ \tilde{i}(r,h):=(i(r),\rho(r)h)$$
$$\tilde{e}(u,h) := (e(u),h) \ \\ \ \ \ \ \& \ \ \ \ \ \\ \tilde{m} \left( (r,h),(r',\rho(r)h) \right) := \left( m(r,r'), \rho(m(r,r'))h \right)$$
where $\rho : \RR \to G$ is a map satisfying:
$$i^* \rho = \rho^{-1} \ \ \ \ \ (\pi_1^* \rho)\cdot (\pi_2^* \rho) = m^* \rho$$
\end{definition}

\begin{definition} For a group $G$ we write $\bar{G}$ to denote
the groupoid $\star \times G \rightrightarrows \star$.
\end{definition}

\begin{proposition}
   To have a principal $G$-bundle over $\GG =(\RR \twoarrows \UU)$ is
   the same thing as to have a morphism of groupoids
   $\GG\to\bar{G}$.
\end{proposition}

This definition coincides with the one of orbibundle given previously in section
\ref{sectionorbifolds} when we work with the groupoid associated to the orbifold,
this will be discussed in detail in the next section.

\section{Orbifolds and groupoids}

\subsection{The groupoid associated to an orbifold} \label{subsectiongroupoidorbifold}

The underlying idea behind what follows is that an orbifold is best understood when
it is interpreted as a stack. We will expand this idea  in the Appendix. There we
explain separately the procedures to go first from an orbifold to a stack, in such a
way that the category of orbifolds constructed above turns out to be a full
subcategory of the category of stacks; and then, from a stack to a groupoid, producing
again an embedding of categories. But there is a more direct way to pass directly
from the orbifold to the groupoid and we explain it now. We recommend to see
\cite{CrainicMoerdijk, Moerdijk2002,  MoerdijkPronk} for a detailed exposition of this issue. Then
we complete the dictionary between the orbifold approach of \cite{ChenRuan} and the
groupoid approach.

 Let $X$ be an orbifold and
$\{(V_p,G_p,\pi_p)\}_{p \in X}$ its orbifold structure, the groupoid $\RR \twoarrows
\UU$ associated to $X$ will be defined as follows: $\UU := \bigsqcup_{p \in X} V_p$
and an element $g: (v_1, V_1) \to (v_2,V_2)$ (an arrow) in $\RR$ with $v_i \in V_i, i
= 1,2$, will be a equivalence class of triples $g= [\lambda_1,w,\lambda_2]$ where $w
\in W$ for another uniformizing system $(W,H,\rho)$, and the $\lambda_i$'s are
injections $(\lambda_i, \phi_i) : (W,H,\rho) \to (V_i,G_i, \pi_i)$ with $\lambda_i(w)
=v_i, i = 1,2$ as in definition \ref{induced}.

For another injection $(\gamma, \psi) : (W',H',\rho') \to (W,H,\rho)$ and $w' \in W'$ with $\gamma(w')=w$ then
$[\lambda_1,w,\lambda_2] = [\lambda_1 \circ \gamma,w',\lambda_2 \circ \gamma]$

Now the maps $s,t,e,i,m$ are naturally described:
$$s([\lambda_1,w,\lambda_2]) = (\lambda_1(w),V_1), \ \ \ \ \  t([\lambda_1,w,\lambda_2]) = (\lambda_2 (w),V_2) \ \ \ \
e(x,V) = [id_V,x,id_V]$$
$$i([\lambda_1,w,\lambda_2]) =[\lambda_2,w,\lambda_1] \ \ \ \ m([[\lambda_1,w,\lambda_2],[\mu_1,z,\mu_2]) =
[\lambda_1 \circ \nu_1 , y, \mu_2 \circ \nu_2]$$
where $h=[\nu_1,y,\nu_2]$ is an arrow joining $w$ and $z$ (i.e. $\nu_1(y) =w \ \& \ \nu_2(y) =z$)

It can  be given a topology to $\RR$ so  that $s,t$ will be \'etale maps, making it
into a proper, \'etale, differentiable groupoid,  and it is not hard to check that all
the properties of groupoid are satisfied \cite[Thm 4.1]{MoerdijkPronk1}.

\begin{rem}
Two equivalent orbifold structures (as in Def. \ref{orbifoldstructure}) 
will induce Morita equivalent groupoids and vice versa. Thus, the choice of groupoid in
the Morita equivalence class
that we will use for a specific orbifold will depend on the setting,
it may change once we take finer covers, but it will be clear that it represents
the same orbifold. 
\end{rem}

This is a good place to note that an orbifold $X$ given by a groupoid
$\RR\rightrightarrows\UU$ will be a smooth manifold if and only if the map
$(s,t)\colon\RR\to\UU\times\UU$ is one-to-one.

Now we can construct principal $\Gamma$-bundles on the groupoid $\RR \twoarrows \UU$
associated to the orbifold $X$ getting,

\begin{proposition}
Principal $\Gamma$ bundles over the groupoid $\RR \twoarrows \UU$ are in 1-1
correspondence with $\Gamma$-orbibundles over $X$.
\end{proposition}

\begin{proof} Let's suppose the bundles are complex, in other
words $\Gamma = \GLC{n}$. The proof for general $\Gamma$ is
exactly the same.

For  an $n$-dimensional complex bundle over $\RR \twoarrows \UU$ we have a map $\rho
: \RR \to \GLC{n}$ and a groupoid structure $\RR \times \complex^n \twoarrows \UU
\times \complex^n$ as in definition \ref{definitionprincipalbundle}. Let $U$ be an
open set of $X$ uniformized by $(V,G, \pi)$ which belongs to its orbifold structure;
for $g \in G$ and $x \in V$, $\xi = [id_G,x,g]$ is an element of $\RR$ ( via the
identity on $V$, and the action of $g$ in $V$ and the conjugation by $g$ on $G$
thought of as an automorphism of $V$) and we can define $\rho_{V,G} : V \times G \to
\GLC{n}$ by $\rho_{V,G} (x,g) \mapsto \rho([id_G,x,g])$. As
$m([id_G,x,g],[id_G,gx,hg])=[id_G,x,hg]$, we have $\rho([id_G,gx,h]) \circ \rho(
[id_G,x,g]) = \rho([id_G,x,hg])$, which implies $\rho_{V,G} (gx,h) \circ
\rho_{V,G}(x,g) = \rho_{V,G} (x,hg)$. So $(V \times \complex^n, G, \tilde{\pi})$ with
$\rho_{V,G}$ extending the action of $G$ in $\complex^n$ is a uniformizing system for
the orbibundle we are constructing, we need to prove now that they define the same
germs and then we would get a orbibundle $E \to X$ using its bundle orbifold
structure.

Let  $(\lambda_i , \phi_i) : (W,H,\mu) \to (V_i,G_i,\pi_i)$ be
injections of uniformizing systems of $X$, with corresponding
bundle uniformizing systems $(W \times \complex^n, H
,\tilde{\mu})$ and $(V_i \times \complex^n, G_i ,
\tilde{\pi_i})$. For $x \in W $, $\xi \in \complex^n$ and
 $h \in H$, $([id_H,x,h] ,\xi) \in \RR \times \complex^n$ and
$\tilde{t}([id_H,x,h] ,\xi) = (hx, \rho([id_H,x,h]) \xi)= \rho_{W,H}(x,h)\xi$. As
 $[\lambda_i,x,\phi_i(h) \circ \lambda_i] = [id_H,x,h]$ for $i \in \{1,2\}$ then
 $\rho_{V_1,G_1}(\lambda_1(x),\phi_1(h)) = \rho_{V_2,G_2}(\lambda_2(x),\phi_2(h))$;
 so the bundle uniformizing
systems $(V_i \times \complex^n, G_i, \tilde{\pi_i})$ define the same germs,
thus they form a bundle orbifold
structure over $X$.

Conversely, if we have the orbibundle  structure for $E \to X$ we need to define the
function $\rho : \RR \to \GLC{n}$. So, for injections $(\tilde{\lambda_i} , \phi_i) :
(W \times \complex^n,H,\mu) \to (V_i \times \complex^n,G_i,\pi_i)$ (where
$\tilde{\lambda_i}$ extends the $\lambda_i$'s previously defined),
$\rho([\lambda_1,x,\lambda_2])$ will be the element in $\GLC{n}$ such that  maps
$pr_2(\tilde{\lambda_1}(x,\xi)) \mapsto pr_2(\tilde{\lambda_2}(x,\xi))$, here $pr_2$
stands for the projection on the second coordinate; in other words
$$ \rho([\lambda_1,x,\lambda_2]) pr_2(\tilde{\lambda_1}(x,\xi)) =pr_2(\tilde{\lambda_2}(x,\xi))$$
Because this bundle uniformizing systems define the same germs,
$\rho$ satisfies the product formula;  the inverse formula is
clearly satisfied.
\end{proof}

\begin{proposition}
Isomorphic $\Gamma$-bundles over $\RR \twoarrows \UU$ correspond to isomorphic $\Gamma$-orbibundles over $X$, and
vice versa.
\end{proposition}
\begin{proof}
We will focus again on complex bundles. To understand what
relevant information we have from isomorphic bundles, let's see
the following lemmas

\begin{lemma} \label{isomorphicbundles}
An isomorphism of bundles  over $\RR \twoarrows \UU $ (with maps $\rho_i : \RR \to \GLC{n}$ for $i=1,2$)
is determined by a map $\delta : \RR \to \GLC{n}$
such that
\begin{displaymath}
\begin{array}{ccc}
\RR \times \complex^n & \stackrel{\Psi}{\To} &\RR \times \complex^n\\
(r, \xi) & \mapsto  & (r, \delta(r)\xi)
\end{array}
\ \ \ \ \ \
\begin{array}{ccc}
\UU \times \complex^n & \stackrel{\psi}{\To} &\UU \times \complex^n\\
(u, \xi) & \mapsto  & (u, \delta(e(u))\xi)
\end{array}
\end{displaymath}
satisfying $\delta(i(r)) \rho_1(r) = \rho_2(r) \delta (r)$ and $\delta(r) = \delta (e  s(r))$.
\end{lemma}
\begin{proof}
It is easy to check that $(\Psi,\psi)$ defined in this way is a morphism between the bundles; the equality
$\delta(r) = \delta (e  s(r))$ comes from the diagram of the \emph{source} map and
 $\delta(i(r)) \rho_1(r) = \rho_2(r) \delta (r)$ from the one of the \emph{target} map, the rest of the diagrams
follow from those two.
\end{proof}
In the same way we could do this procedure for complex orbibundles:
\begin{lemma}
An isomorphism of complex orbibundles over $X$ (with maps $\rho^i_{V,G} :V \times G \to \GLC{n}$ for $i=1,2$ and
$\{(V,G,\pi)\}$ orbifold structure of $X$) is determined by the maps $\tilde{\delta}_V:V \to \GLC{n}$ such that
\begin{displaymath}
\begin{array}{ccc}
V \times \complex^n & \to & V \times \complex^n \\
 (r,\xi) & \mapsto & (r,\tilde{\delta}_V(r)\xi)
\end{array}
\end{displaymath}
satisfying $\tilde{\delta}(gr) \rho^1_{V,G}(r,g) = \rho^2_{V.G}(r,g) \tilde{\delta} (r)$. The $\tilde{\delta}_V$'s
 form a good map.
\end{lemma}
\begin{proof}
Because the underlying orbifold structure needs to be mapped to itself, we obtain the $\delta_V$'s.
The equality $\tilde{\delta}(gr) \rho^1_{V,G}(r,g) = \rho^2_{V.G}(r,g) \tilde{\delta} (r)$ holds
 because the good map condition.
\end{proof}
The proof of the proposition is straight forward from these lemmas. The map $\delta$
that comes from the isomorphism of the complex bundles determines uniquely the
$\tilde{\delta}_V$'s, and vice versa.
\end{proof}

\begin{example}\label{framebundle}
   The \emph{tangent bundle} $TX$ of an orbifold $X$ is a
   orbibundle over $X$. If $U=V/G$ is a local uniformizing system, then a
   corresponding local uniformizing system for $TX$ will be $TU/G$ with
   the action $g\cdot(x,v)=(gx,dg_x(v))$.

   Similarly the \emph{frame bundle} $P(X)$ is a principal orbibundle over
   $X$. The local uniformizing system is $U\times\GLC{n}/G$ with local
   action $g\cdot(x,A)=(gx,dg\circ A)$. Notice that $P(X)$ is
   always a \emph{smooth manifold} for the local action is free
   and $(s,t)\colon\RR\to\UU\times\UU$ is one-to-one.
\end{example}

We  want the morphism between orbifolds to be morphisms of groupoids, and this is
precisely the case for the good maps given in Definition \ref{definitiongoodmap}

\begin{proposition}
 \label{morphism=goodmap}
A morphism of groupoids induces a good map between the underlying orbifolds, and
conversely, every good map arises in this way.
\end{proposition}

\begin{proof}
For ${f} : X \to X'$  a good
map between orbifolds, we have a correspondence $U \leftrightarrow U'$ between open subsets of a
compatible cover of $X$ and  open subsets of $X'$, such that $f(U) \subset U'$ , and
$U_1 \subset U_2$ implies $U_1' \subset U_2'$. Moreover, we are provided with local
liftings $f_{UU'}:(V,G, \pi) \to (V', G', \pi')$ as in the Definition
\ref{definitiongoodmap}.
Let $\RR \twoarrows \UU$ and $\RR' \twoarrows \UU'$ be the groupoids 
constructed from the orbifold structures of $X$ and $X'$
respectively, determined by the compatible cover $\{U_i\}_i$ of $X$ and a cover of $X'$
that uniformizes $\{{U'_j}\}_j$.

Define  $\psi : \UU \to \UU'$ such that $\psi |_U = f_{UU'}$ and  $\Psi :\RR \to \RR'$
by  $\Psi([\lambda_1, w ,\lambda_2]) = ([\nu(\lambda_1), \psi(w)
, \nu(\lambda_2)]$, where the $\lambda_i$'s are injections
between $W$ and $V_i$ and the $\nu(\lambda_i)$'s are the
corresponding injections between $W'$ and $V_i'$ given by the
definition of good map; because
 $$ \nu(\lambda_i) \circ f_{WW'} = f_{V_iV_i'} \circ \lambda_i$$
the  function $\Psi$ is well defined and together with $\psi$, satisfy all the
conditions for a morphism of groupoids.

It is clear that the groupoids just used could differ from the groupoids one obtain
after performing the construction defined at the beginning of this chapter, but they
are respectively Morita equivalent. 

On the  other hand, if we are given $\Psi: \RR \to \RR'$ and
$\psi : \UU \to \UU'$, we can take a sufficiently small open
compatible cover for $X$ such that for $U$ in its cover there is
an open set $U'$ of $X'$ with the desired properties. For
$(V,G,\pi)$ and $(V',G',\pi')$ uniformizing systems of $U$ and
$U'$ respectively, we need to define $f_{UU'}$. The map between
$V$ and $V'$ is given by $\psi |_V$, and the injection between
$G$ and $G'$ is given as follows.

Let's  take $x \in V$ and $g \in G$; we have an automorphism of $(V,G,\pi)$ given by
the action by $g$ on $V$ and by conjugation on $G$, call this automorphism
$\lambda_g$; then $[id_G, x, g]$ is an element of $\RR$, using the properties of
$\Psi$ and $\psi$ we get that $\Psi([Id_G,x, g])=[Id_{G'},\psi(x),g']$, where $g' \in
G'$; this because every automorphism of $(V',G',\pi')$ comes from the action of an
element in $G'$ (see \cite[Lemma 2.1.1]{Ruan}); moreover,  we have that $g' \circ
\psi(x) = \psi \circ g(x)$. This will give us an homomorphism $\rho_{UU'}:G \to G'$
sending $g \mapsto g'$ that together with $f_{UU'}$ form the compatible system we
required.
\end{proof}

Orbifolds have the property that they can be seen as the quotient of a
manifold by a Lie group. We just construct the frame bundle $P(X)$ of $X$, which is a
manifold, together with the natural action of $\Or{n}$ as in Example
\ref{framebundle} (cf. \cite{AdemRuan}).
\begin{example} \label{exampleorbifolduniversalcover}

Let $X$  be a $n$-dimensional orbifold, $Y$ its orbifold universal cover and
$H=\pi_1^{orb}(X)$ its fundamental orbifold group and $f\colon Y \to X$ the cover
good map. Let $P(Y)$ be the frame bundle of $X$, by \ref{framebundle} we know that
$P(X)$ is a smooth manifold and it is endowed with a smooth and effective $\Or{n}$
action with finite isotropy subgroups such that $X\simeq [P(X)/\Or{n}]$ in the
category of orbifolds (cf. \cite{AdemRuan} prop. 2.3.)

The frame bundle  $P(Y)$ is isomorphic to $f^*P(X)$ and lifting the action of $H$ in
$Y$ to a free action of $H$ in $P(Y)$ with $P(Y)/H\simeq P(X)$ we obtain the
following diagram.
        $$\xymatrix{
        P(Y) \ar[r]^{f'}_{/H} \ar[d]_{/O(n)} & P(X) \ar[d]^{/O(n)} \\
        Y \ar[r]^f_{/H} & X}
        $$

Let's consider now the groupoids $\RR_Y \stackrel{s_Y,t_Y}{\twoarrows} \UU_Y$ and
$\RR_X \stackrel{s_X,t_X}{\twoarrows} \UU_X$ associated to the orbifolds $Y$ and $X$
by using their frame bundles (i.e. $\RR_Y = P(Y) \times \Or{n}$ and $\UU_Y=P(Y)$ with
$s_Y$ and $t_Y$ as in Example \ref{quotient}.) Every $h\in H$ induces a morphism
of groupoids
        $$\xymatrix{
        \RR_Y \ar[r]^h \ar@<-.5ex>[d] \ar@<+.5ex>[d] & \RR_X \ar@<-.5ex>[d] \ar@<+.5ex>[d]\\
        \UU_Y \ar[r]^h & \UU_X}
        $$

since the action of $h$ in $P(Y)\times \Or{n}$ commutes with the action of $\Or{n}$,
for $P(Y)$ is simply $f^*P(X)$.

As we are  working in the reduced case, the orbifold structures of $Y$ and $X$ can be
obtained using the frame bundles $P(Y)$ and $P(X)$ so we can choose a sufficiently
small orbifold cover $\{U\}$ of $Y$, such that for $(V,G,\pi)$ a uniformizing system
of $U$, and $h \in H$ ,we have an isomorphism $\eta_h: (V,G, \pi) \cong
(V',G,'\pi')$, where $(V',G', \pi')$ is a uniformizing system for $U'=hU$. In other
words, the map $\eta_h$ induces a groupoid automorphism of the orbifold (a good map).

Let $\RR_Y \times H \stackrel{s,t}{\twoarrows} \UU_Y$ be the groupoid defined by the
following maps
$$s(r,h)= s_Y(r) \ \ \ \ \ t(r,h) = h(t_Y(r)) \ \ \ \ \ \ e(x) = (e_Y(x),id_H)$$
$$i(r,h) = (h(i_Y(r)),h^{-1}) \ \ \ \ \ m((r_1,h_1),(r_2,h_2))= (m(r_1,h^{-1}(r_2)),h_2 h_1)$$
then the following holds.

\begin{proposition}
The groupoids $\RR_Y \times H \twoarrows \UU_Y$ and $\RR_X \twoarrows \UU_X$ are
Morita equivalent.
\end{proposition}

\begin{proof}
Noting that the map $P(Y)\to P(X)$ is a surjection and recalling that $\RR_X = P(X)
\times \Or{n} $ we can see that $s \circ \pi_2 \colon \UU_Y {}_f \times_t \RR_X \to
\UU_X$ is an \'etale surjection. Finally because the action of $H$ in $\RR_Y$ and
$\UU_Y$ is free and $\RR_Y / H \simeq \RR_X$, $\UU_Y/ H \simeq \UU_X$ it is immediate
to verify that
        $$\xymatrix{
        \RR_Y \times H \ar[r]^{f' \pi_1} \ar[d]_{(s,t)} & \RR_X \ar[d]^{(s_X,t_Y)}\\
        \UU_Y \times \UU_Y \ar[r]^h & \UU_X \times \UU_X }
        $$
is a fibered square.
\end{proof}

\end{example}

\subsection{The category associated to a groupoid and its classifying space}

To every groupoid $\RR \twoarrows \UU$ we can  associate a category $\CC$ whose
objects are the objects in $\UU$ and whose morphisms are the objects in $\RR$ that we
have called arrows before. We can see
$$\RR^{(n)} := \underbrace{\RR \timests \cdots \timests
\RR}_{n}$$ as the composition of $n$ morphisms. In the case in
which $\RR$ is a set then $\RR^{(n)}$ is the set of sequences
$(\gamma_1, \gamma_2, \ldots, \gamma_n)$ so that we can form the
composition $\gamma_1 \circ \gamma_2 \circ \cdots \circ \gamma_n$.

We this data we can form a simplicial set \cite{MaySimplicial,
Segal1}.

\begin{definition} A \emph{semi-simplicial set (resp. group, space, scheme)}
$X_\bullet$ is a sequence of sets $\{X_n\}_{n\in \naturals}$
(resp. groups, spaces, schemes) together with maps
\begin{equation*}
    X_0 \leftrightarrows X_1 \leftrightarrows X_2
    \leftrightarrows \cdots \leftrightarrows X_m \leftrightarrows
    \cdots
\end{equation*}
\begin{equation}
   \partial_i\colon X_m \to X_{m-1}, \ \ \ \ s_j\colon X_m\to
   X_{m+1}, \ \ \ \ 0\leq i,j \leq m.
\end{equation}
called \emph{boundary} and \emph{degeneracy} maps, satisfying
\begin{eqnarray*}
\partial_i \partial_j &=& \partial_{j-1} \partial_{i} \ \ \ \mbox{if $i<j$} \\
 s_i s_j &=& s_{j+1} s_i \ \ \ \mbox{if $i<j$} \\
\partial_i s_j &=& \left\{ \begin{array}{ll}
                            s_{j-1} \partial_i & \mbox{if $i<j$}\\
                            1 & \mbox{if $i=j,j+1$}\\
                            s_j \partial_{i-1} & \mbox{if
                            $i>j+1$}\\
                           \end{array}
                   \right.
\end{eqnarray*}
\end{definition}

The nerve of a category (see \cite{Segal1}) is a semi-simplicial set $N\CC$ where the
objects of $\CC$ are the vertices, the morphisms the 1-simplices, the triangular
commutative diagrams the 2-simplices, and so on. For a category coming from a
groupoid then the corresponding simplicial object will satisfy
$N\CC_n=X_n=\RR^{(n)}$. We can define the boundary maps $\partial_i : \RR^{(n)} \to
\RR^{(n-1)}$ by:
\begin{displaymath}
\partial_i(\gamma_1, \dots , \gamma_n) = \left\{
 \begin{array}{ll}
(\gamma_2, \dots , \gamma_n) & \mbox{if $ i=0$} \\
(\gamma_1, \dots, m(\gamma_i,\gamma_{i+1}), \dots , \gamma_n) & \mbox{ if $ 1 \leq i \leq n-1$} \\
(\gamma_1, \dots, \gamma_{n-1}) & \mbox{if $ i =n $} \\
\end{array}
\right.
\end{displaymath}
and the degeneracy maps by
\begin{equation*}
   s_j(\gamma_1,\ldots,\gamma_n)=\left\{ \begin{array}{ll}
                                          (e(s(\gamma_1)),\gamma_1,\ldots,\gamma_n)&\mbox{for
                                          $j=0$} \\
                                          (\gamma_1,\ldots,\gamma_j,e(t(\gamma_j)),
                                          \gamma_{j+1},\ldots,\gamma_n)&\mbox{for
                                          $j\geq1$}
                                         \end{array}
                                  \right.
\end{equation*}

We will write $\Delta^n$ to denote the standard $n$-simplex in
$\real^n$. Let $\delta_i\colon\Delta^{n-1}\to\Delta^n$ be the
linear embedding of $\Delta^{n-1}$ into $\Delta^n$ as the $i$-th
face, and let $\sigma_j\colon\Delta^{n+1}\to\Delta^{n}$ be the
linear projection of $\Delta^{n+1}$ onto its $j$-th face.

\begin{definition} The \emph{geometric realization} $|X_\bullet|$
of the simplicial object $X_\bullet$ is the space
\begin{equation*}
|X_\bullet|=\left.\left(\coprod_{n\in\naturals} \Delta^n \times
X_n \right)\right/ \begin{array}{c}
                    (z,\partial_i(x))\sim(\delta_i(z),x)\\
                    (z,s_j(x))\sim(\sigma_j(z),x)
                   \end{array}
\end{equation*}
Notice that the topologies of $X_n$ are relevant to this
definition.

The semi-simplicial object $N \CC$ determines $\CC$ and its
topological realization is called $B \CC$, the \emph{classifying
space of the category}. Again in our case $\CC$ is a
\emph{topological category} in Segal's sense \cite{Segal1}.

For a groupoid $\RR\rightrightarrows\UU$ we will call
$B(\RR\rightrightarrows\UU)=B\CC=|N\CC|$ the \emph{classifying space of the groupoid}.

\end{definition}

The following proposition establishes that $B$ is a functor from
the category of groupoids to that of topological spaces. Recall
that we say that two morphisms of groupoids are Morita related if
the corresponding functors for the associated categories are
connected by a morphism of functors.

\begin{proposition}[\cite{Moerdijk98}] (cf. \cite[Prop. 2.1]{Segal1}) \label{independence1}
   A morphism of groupoids $\Xx_1 \to \Xx_2$ induces a continuous
   map $B\Xx_1 \to B\Xx_2$. Two morphism that are Morita
   related will produce homotopic maps. In particular a Morita
   equivalence $\Xx_1 \sim \Xx_2$ will induce a homotopy
   equivalence $B\Xx_1\simeq B\Xx_2$. This assignment is
   functorial.
\end{proposition}

\begin{example} For the groupoid $\bar{G}=(\star \times G
\rightrightarrows \star)$ the space $B\bar{G}$ coincides with the
classifying space $BG$ of $G$.

Consider now the groupoid $\Xx=(G\times G \rightrightarrows G)$ where $s(g_1,
g_2)=g_1$, $t(g_1,g_2)=g_2$ and $m((g_1,g_2);(g_2,g_3))=(g_1,g_3)$ then it is easy to
see that $B\Xx$ is contractible and has a $G$ action. Usually $B\Xx$ is written $EG$

A morphism of groupoids $\Xx \to \bar{G}$ is the same thing as a principal $G$ bundle
over $\Xx$ and therefore can be written by means of a map $G \times G \to G$. If we
choose $(g_2,g_2)\mapsto g_1^{-1} g_2$ the induced map of classifying spaces
$$ EG \longrightarrow BG $$
is the universal principal $G$-bundle fibration over $BG$.
\end{example}

\begin{example}
Consider a smooth manifold $X$ and a good open cover $\UU=\{U_{\alpha}\}_\alpha$.
Consider the groupoid $\GG=(\RR \rightrightarrows \UU)$ where $\RR$ consists on the
disjoint union of the double intersections $U_{\alpha\beta}$. Segal \cite[Prop.
4.1]{Segal1} calls $X_U$ the corresponding topological category. There he proves that
$B\GG = BX_U \simeq X$.

If we are given a principal $G$ bundle over $\GG$ then we have a morphism
$\GG\to\bar{G}$ of groupoids, that in turn induces a map $X\to BG$. Suppose that in
the previous example we take $G=\GLC{n}$. Then we get a map $X\to B\GLC{n}=BU$ and
since $K(X)=[X,BU]$ this is an element in $K$-theory.
\end{example}

\begin{example}
   Consider a groupoid $\Xx$ of the form $M\times G \rightrightarrows M$ where $G$ is
   acting on $M$ continuously. Then $B\Xx \simeq EG\times_G M$ is the Borel
   construction for the action $M \times G \to M$.
\end{example}

\subsection{Sheaf Cohomology and Deligne Cohomology}\label{sheaf}

On a smooth manifold $X$ a sheaf $S$ can be defined as a functor from the category
whose objects are open sets of $X$, and whose morphisms are inclusions to the category
(for example) of abelian groups, and a gluing condition of the type described in the
Appendix. So for every open set $U_\alpha$ in $X$ we have an abelian group
$S_{U_\alpha}=\Gamma_S(U_\alpha)$ called the sections of $S$ in $U_\alpha$. In the
representation of a smooth manifold as a groupoid $\RR \rightrightarrows \UU$ where
$$\UU = \bigsqcup_\alpha U_\alpha \ \ \  \ \RR= \bigsqcup_{(\alpha, \beta)} U_\alpha \cap
U_\beta \ \ (\alpha,\beta) \neq (\beta,\alpha).$$  A sheaf can be encoded by giving a
sheaf over $\UU$ with additional gluing conditions given by $\RR$.

\begin{definition}\cite{CrainicMoerdijk}
A \emph{sheaf} $\SSS$ on a groupoid $\RR \rightrightarrows \UU$ consists of
\begin{enumerate}
\item A sheaf $S$ on $\UU$.
\item A continuous (right) action of $\RR$ on the total space of $S$.
\end{enumerate}
An \emph{action} of $\RR$ on $S\toparrow{\pi} \UU$ is a map $S {\: {}_{\pi}  \!
\times_{t}} \RR \to S$ satisfying the obvious identities.
\end{definition}

The theory of sheaves over groupoids and their cohomology has been developed by
Crainic and Moerdijk \cite{CrainicMoerdijk}. There is a canonical notion of morphism
of sheaves. So we can define the category $\SSS{h}(\Xx)$ of sheaves over the groupoid
$\Xx$. Morita equivalent groupoids have equivalent categories of sheaves. There is a
notion of sheaf cohomology of sheaves over groupoids defined in terms of resolutions.
There is also a \v{C}ech version of this cohomology developed by Moerdijk and Pronk
\cite{MoerdijkPronk}. We will call a groupoid with a sheaf a \emph{Leray groupoid} if
the associated sheaf cohomology can be calculated as the corresponding \v{C}ech
cohomology. From now on we will always take a representative of the Morita class of
the groupoid that is of Leray type. The basic idea is just as in the case of a smooth
manifold, an $S$ valued $n$ \v{C}ech cocycle is an element $c \in \Gamma_S(\coprod
U_{\alpha_1 \cdots \alpha_n})$, and in a similar fashion we can define cocycles in
the groupoid $\RR \rightrightarrows \UU$ in terms of the sheaf and the products
$\RR^{(n)} := \underbrace{\RR \timests \cdots \timests \RR}_{n}$. Then using
alternating sums of the natural collection of maps
\begin{displaymath}
   \RR^{(0)} \leftleftarrows \RR^{(1)} \overleftarrow{\leftleftarrows} \RR^{(2)}
   \underleftarrow{\overleftarrow{\leftleftarrows}} \RR^{(3)} \cdots
\end{displaymath}
we can produce boundary homomorphisms and define the cohomology theory.

The resulting groupoid sheaf cohomology satisfies the usual long exact sequences and
spectral sequences. In particular we can use the exponential sequence induced by the
sequence of sheaves $0\to\integer\toparrow{i}\real\toparrow{\exp}\U{1}\to1$.

In \cite{MoerdijkPronk, Moerdijk, CrainicMoerdijk} we find a
theorem that implies the following

\begin{theorem}\label{Pronk} For an orbifold with groupoid $\Xx$ and a
locally constant system $A$ of coefficients (for example
$A=\integer$) we have
$$ H^*(\Xx, A ) \cong H^*(B\Xx, A )$$
where the left hand side is orbifold sheaf cohomology and the right hand side is
ordinary simplicial cohomology.
\end{theorem}

Moerdijk has proved that the previous theorem is true for
arbitrary coefficients $A$ \cite{Moerdijk98}.

Crainic and Moerdijk have also defined hypercohomology for a bounded complex of
sheaves in a groupoid, and they obtained the basic spectral sequence. 
 In \cite{LupercioUribe4, LupercioUribe3} we define \emph{Deligne cohomology} for groupoids
associated to orbifolds and also \emph{Cheeger-Simons cohomology}.
\section{Gerbes over orbifolds}

\subsection{Gerbes and inner local systems.}

From this section on we are going to work  over the groupoid associated to an
orbifold. For $\RR \twoarrows \UU$ the groupoid associated to an orbifold $X$ defined
in \ref{subsectiongroupoidorbifold} we will consider the following

\begin{definition}\label{gerbedef}
A gerbe over an orbifold $\RR \twoarrows \UU$,  is a complex line
bundle $\LL$ over $\RR$ satisfying the following conditions
\begin{itemize}
\item $i^* \LL \cong \LL^{-1}$
\item $\pi_1^* \LL \otimes \pi_2^* \LL \otimes m^*i^* \LL \stackrel{\theta}{\cong} 1$
\item $\theta : \RR \timests \RR \to \U{1}$ is a 2-cocycle
\end{itemize}
where $\pi_1,\pi_2 : \RR \timests \RR \to \RR$ are  the
projections on the first and the second coordinates, and $\theta$
is a trivialization of the line bundle.
\end{definition}

The following proposition\footnote{We owe this observation to I.
Moerdijk} shows that the analogy with a finite group can be
carried through in this case.

\begin{proposition}\label{gerbeextension} To have a gerbe $\LL$ over a groupoid $\GG$ is
the same thing as to have a central extension of groupoids
$$ 1 \to \overline{\U{1}} \to \tilde{\GG} \to \GG \to 1 $$
\end{proposition}

\begin{lemma}
In the case of a smooth manifold $M$ \cite{Witten, Hitchin} we define the groupoid as
in the example \ref{examplemanifold}, for a line bundle $\LL$ over $\RR$ we get line
bundles $\LL_{\alpha \beta} := \LL|_{U_{\alpha \beta}}$ over the double intersections
$U_{\alpha \beta}$ such that $\LL_{\alpha \beta} \cong \LL_{\beta \alpha}^{-1}$, and
$\LL_{\alpha \beta} \LL_{\beta \gamma} \LL_{\alpha \gamma}^{-1}
\stackrel{\theta}{\cong} 1$ over the triple intersections $U_{\alpha \beta \gamma}$;
then we get a gerbe over the manifold as defined in section
\ref{sectiongerbesmanifolds}.
\end{lemma}

We want to relate the discrete torsions of Y. Ruan \cite{Ruan} over a discrete group
$G$ and the gerbes over the corresponding groupoid

\begin{example} \label{gerbesdiscretegroup}
Gerbes over a discrete group $G$ are in 1-1 correspondence  with
the set of two-cocycles $Z(G, \U{1})$.

 We recall that $\overline{G}$ denotes the groupoid $* \times G \twoarrows *$ the
trivial maps $s,t$ and $i(g) = g^{-1}$ and $m(h,g) = hg$ (clearly we can drop the $*$
as it is customary). A gerbe over $\overline{G}$ is a line bundle $\LL$ over $G$ such
that, if we call $\LL_g$ the fiber at $g$, $\LL_g^{-1} = \LL_{g^{-1}}$ and $\LL_g
\LL_h \stackrel{\beta}{\cong} \LL_{gh}$. So for each $g,h \in G$ we have a
trivialization $\beta_{g,h} \in \U{1}$ satisfying
$\beta_{g.h}\beta_{gh,k}=\beta_{g,hk}\beta_{h,k}$ because
$$\LL_g \LL_h \LL_k \stackrel{\beta_{g,h}}{\cong} \LL_{gh}\LL_k \stackrel{\beta_{gh,k}}{\cong} \LL_{ghk}
\stackrel{\beta_{g,hk}^{-1}}{\cong} \LL_g \LL_{hk} \stackrel{\beta_{h,k}^{-1}}{\cong}
\LL_g \LL_h \LL_k$$ Then $\beta : G \times G \to \U{1}$ satisfies the cocycle
condition and henceforth is a two-cocycle.

It is clear how to construct the gerbe over $G$ once we have the two-cocycle.
\end{example}

The representations of $L^\alpha_g:C(g) \to \U{1}$ defined  in
section \ref{sectionorbifoldcohomology} for some $\alpha \in
H^2(G,\U{1})$ come from the fact that
$$\LL_g \LL_h  \stackrel{\alpha_{g,h} \alpha^{-1}_{h,g}}{\cong} \LL_{h}\LL_{g}; $$
 then $\theta(g,h) := \alpha_{g,h} \alpha_{h,g}^{-1}$ defines a
representation $\theta(g, \_): C(g) \to \U{1}$ and it matches the $L_g^\alpha$ for
$\beta = \alpha$.

\begin{rem}
Every inner local system over an orbifold $X$ defined by Ruan \cite{Ruan}  as in
Definition \ref{innerlocalsystem}, comes from a gerbe on the groupoid $\RR \twoarrows
\UU$ associated to it. This is because $\RR$ contains copies of the twisted sectors.
This is explained in detail in \cite{LupercioUribe2}
\end{rem}

\subsection{The characteristic class of a Gerbe}

   We want to classify gerbes over an orbifold.
   As we have pointed out before the family of isomorphism classes
of gerbes on a groupoid $\Gpd$ forms a group under the operation of tensor product of
gerbes, that we will denote as $\Gerbe(\Gpd)$. Given an element $[\LL]\in\Gerbe(\Gpd)$
we can choose a representative $\LL$ and such representative will have an associated
cocycle $\theta \colon \RR \timests \RR \to U(1)$. Two isomorphic gerbes will differ
by the co-boundary of a cocycle $\RR \to U(1)$

\begin{example}
$\Gerbe(\bar G) \cong H^2(G,\U{1})$

Using  lemma \ref{gerbesdiscretegroup} and the previous definition of the group
$\Gerbe(\bar G)$ we see that two isomorphic gerbes define cohomologous cycles, and
vice versa.
\end{example}

We will call the cohomology class $\langle\LL\rangle \in
H^2(\Gpd,\underline{\complex^*})$ of $\theta$, the \emph{characteristic class} of the
gerbe $\LL$. As explained in Section \ref{sheaf} we can use the exponential sequence
of sheaves to show that $H^2(\Gpd,\underline{\complex^*})\cong H^3(\Gpd,\integer)$
and then using the isomorphism \ref{Pronk} $H^3(\Gpd,\integer) \cong
H^3(B(\Gpd),\integer)$ we get

\begin{proposition} \label{character} For a groupoid $\Gpd$ we have the following isomorphism
$$ \Gerbe(\Gpd) \cong H^3(B(\Gpd),\integer)$$
given by the map $[\LL]\to \langle\LL\rangle$ that associates to a gerbe its
characteristic class.
\end{proposition}

In particular using \ref{independence1} we have that

\begin{proposition} \label{indepchar}
   The group $ \Gerbe(\Gpd) $ is independent of the Morita class of $\Gpd$.
\end{proposition}

This could also have been obtained noting that a gerbe over an orbifold can be given
as a sheaf of groupoids in the manner of \ref{gerbesassheaves}

\begin{example}
   Consider an inclusion of (compact Lie) groups $K\subset G$ and consider the groupoid
   $\GG$ given by the action of $G$ in $G/K$,
   $$ G/K \times G \rightrightarrows G/K $$

   Observe that the stabilizer of $[1]$ is $K$ and therefore we have that the
   following groupoid
   $$ [1]\times K \rightrightarrows [1] $$
   is Morita equivalent to the one above.

   From this we obtain
   $$ \Gerbe(\GG) \cong H^3(K,\integer) $$
\end{example}

As it was explained in \ref{smoothtwist} in the case of a smooth manifold $X$ we have
that
$$ \Gerbe(X) = [ X , BB\complex^* ] $$
where $BB\complex^* = B\Proj U(\HH)$ for a Hilbert space $\HH$.

Let us write $\overline{\PP U(\HH)}$ to denote the groupoid $\star \times \Proj U(\HH) \to
\star$. We have the following

\begin{proposition}\label{GerbeAsHomo}
   For an orbifold $X$ given by a groupoid $\Xx$ we have
   $$ \Gerbe(\Xx) = [\Xx , \overline{\PP U(\HH)}] $$
   where $[\Xx , \overline{\PP U(\HH)}]$ represents the Morita equivalence classes of morphisms
   from $\Xx$ to $\overline{\PP U(\HH)}$
\end{proposition}

\subsection{Differential geometry of Gerbes over orbifolds and the $B$-field.}

Just as in the case of a gerbe over a smooth manifold, we can do differential geometry
on gerbes over an orbifold groupoid $\Xx=(\Gpd)$. Let us define a connection over a
gerbe in this context.

\begin{definition}
   A \emph{connection} $(g,A,F,G)$ over a gerbe $\Gpd$ consists of a complex valued
   0-form $g\in\Omega^0 (\RR
   \timests \RR)$, a 1-form $A\in \Omega^1(\RR)$, a 2-form $F\in\Omega^2(\UU)$ and a
   3-form $G\in\Omega^3(\UU)$ satisfying
   \begin{itemize}
      \item $G=dF$,
      \item $t^* F - s^* F = dA$ and
      \item $\pi_1^* A + \pi_2^* A +  m^* i^* A = -\sqrt{-1} g^{-1} dg$
   \end{itemize}
   The 3-form $G$ is called the $curvature$ of the connection. A connection is called
   $flat$ if its curvature $G$ vanishes.
\end{definition}

The 3-curvature $\frac{1}{2\pi\sqrt{-1}}G$ represents the integer characteristic class
of the gerbe in cohomology with real coefficients, this is the Chern-Weil theory for a
gerbe over an orbifold. One can reproduce now Hitchin's arguments in \cite{Hitchin}
\emph{mutatis mutandis}. In particular when a connection is flat one can
speak of a holonomy class in $H^2(B\Xx,U(1))$. Hitchin's discussion relating a gerbe
to a line bundle on the loop space has an analogue that we have studied in
\cite{LupercioUribe2}. There, for a given groupoid $\Xx$ we construct a groupoid
$\Loop \Xx$ that represents the free loops on $\Xx$. The `coarse moduli space' or
quotient space of this groupoid coincides with Chen's definition of the loop space
\cite{ChenLoop}, but $\Loop \Xx$ has more structure. In particular if we are given a
gerbe $L$ over $\Xx$, using the holonomy we construct a 'line bundle' $\Lambda$ over
$\Loop \Xx^{S^1}$, the fixed subgroupoid under the action of $S^1$, by a groupoid
homomorphism $\Loop\Xx^{S^1} \To \overline{\U{1}}$. Let us consider the
 groupoid $\wedge \Xx= (\wedge\Xx)_1 \rightrightarrows(\wedge\Xx)_0$,
with objects $(\wedge\Xx)_0=\{r\in\RR | s(r)=t(r)\}$, and arrows
$(\wedge\Xx)_1=\{\lambda\in\RR | r_1\toparrow{\lambda} r_2
\Leftrightarrow m(\lambda,r_2) = m(r_1, \lambda)\}$. The groupoid
$\wedge \Xx$ is certainly \'etale, but it is not necessarily
smooth. In other words the twisted sectors are an \emph{orbispace}
or a topological groupoid \cite{ChenLoop, LupercioUribe2}.
\begin{theorem}\cite{LupercioUribe2}\label{holonomy} The orbifold
 $\widetilde{\Sigma_1 X}$ defined in
\ref{sectionorbifoldcohomology} is represented by the groupoid
$\wedge \Xx$. There is a natural action of $S^1$ on $\Loop\Xx$.
The fixed subgroupoid $(\Loop\Xx)^{S^1}$ under this action is
equal to $\wedge \Xx$. The holonomy line bundle
$\Lambda$ over $\wedge \Xx$ is an inner local system as defined in
\ref{sectionorbifoldcohomology}.
\end{theorem}

From this discussion we see that in orbifolds with discrete torsion as the ones
considered by Witten in \cite[p. 34]{Witten}, what corresponds to the $B$-field 3-form
$H$ in \cite[p. 30]{Witten} is the 3-form $G$ of this section. The analogue of
$K_{[H]}$ that Witten requires in \cite[p. 34]{Witten} will be constructed in the
next section.

Let us recall that the smooth Deligne cohomology groups of an orbifold $\Xx$ can be
defined as in section \ref{sheaf}. To finish this section let us state one last
proposition in the orbifold case.

\begin{proposition} \cite[Prop. 3.0.6]{LupercioUribe3}
   The group of gerbes with connection over an orbifold $\Xx$ are classified by the
   Deligne cohomology group $H^3(\Xx,\integer(3)^\infty_D)$.
\end{proposition}

\section{Twisted $\Kgpd{\LL}$-theory}

\subsection{Motivation}

Think of a group $G$ as the groupoid $\star \times G
\rightrightarrows \star$, that is to say, a category with one
object $\star$ and an arrow from this object to itself for every
element $g\in G$:
$$
 \xymatrix@ -5pt{
 \star  \ar@(r,u)[]|{g}
 }
$$

 The theory of representations of $G$ consists
of the study of the functor $G\mapsto R(G)$ where $R(G)$ is the
Grothendieck ring of representations of $G$ with direct sum and
tensor product as operations.

An $n$-dimensional representation $\rho$ of $G$ is a continuous
assignment of a linear map
$$
 \xymatrix@ -5pt{
 \complex^n  \ar@(r,u)[]|{\rho_g}
 }
$$
 for every arrow $g \in G$. Namely a representation is
encoded in a map
$$ \rho \colon \RR = \star \times G \to \GLC{n}$$
or in other words is a principal $\GLC{n}$ bundle over the groupoid
$\overline{G}=(\star \times G \rightrightarrows \star)$. For finite groups this is
simply an orbibundle.

For a given orbifold, the study of its $\Korb{}$-theory is the exact analogue to the
previous situation, in other words, the representation theory of groupoids is
$K$-theory. The analogue of a representation is an orbibundle as in
\ref{definitionprincipalbundle}. Every arrow in the groupoid corresponds to an
element  in $\GLC{n}$ but now there are many objects so we get a copy of $\complex^n$
for every object in $\UU$, namely a bundle over $\UU$ with gluing information.

In the case of a smooth manifold this recovers the usual
$K$-theory.

It is clear now that we can twist $\Kgpd{}(X)$ by a gerbe $\LL$ over $X$ in the very
same manner in which $R(G)$ can be twisted by an extension
$$ 1\to\complex^*\to\tilde{G}\to G \to 1$$
such an extension is the same thing that a gerbe over $\overline{G} = (\star \times G
\rightrightarrows \star)$. This twisting recovers all the twistings of $K$-theory
mentioned before in this paper.


For a moment lets restrict our attention to the groupoid $\GG$
associated to a smooth manifold $X$ defined in
\ref{examplemanifold}, and let's see how its $K$-theory can be
interpreted in terms of this groupoid.

Let $\CC$ be the (discrete) category whose objects are finite
dimensional vector spaces and whose morphisms are linear
mappings. Then a functor of categories
$$ \GG \longrightarrow \CC $$
assigns to every object of $\GG$ a vector space and to every
morphism of $\GG$ a linear isomorphism in a continuous fashion.
If we recall that the groupoid $\GG$ is given by
$$\UU = \bigsqcup_\alpha U_\alpha \ \ \
  \ \RR= \bigsqcup_{(\alpha, \beta)} U_\alpha \cap U_\beta \ \ (\alpha,\beta) \neq (\beta,\alpha)$$
then we realize that this is equivalent to giving a trivial
vector over $\UU$ and linear gluing instructions, that is to say,
a vector bundle over $M$.

It is also clear that the category $\CC$ is equivalent to the
category with one object for every non-zero integer $n \in
\integer_{\geq 0}$, and with morphisms generated by the
isomorphisms $\bigsqcup \GLC{n}$ and the arrows $n\to m$ whenever
$n \leq m $. The classifying space of $\CC$ is
$\Grr(\complex^\infty)\simeq B\GLC{\infty}\simeq BU$.

In this case $B\GG\simeq X$ \cite{Segal1}, and we get an element
in the reduced $K$-theory of $X$, $[X,BU]=\widetilde{K}(X)$. This
discussion is also valid in the case in which $X$ is an orbifold
and shows that our constructions do not depend on the choice of
Leray \'etale groupoid representing the orbifold that we take.
This will motivate us \emph{define} the $K$-theory of an orbifold
$X$ given by a groupoid $\GG$ by means of such functors
$\GG\to\CC$. Actually we will use groupoid homomorphisms from
$\GG$ to some groupoid $\VV$ whose classifying space is homotopic
to $\CC$, and this will allow us to generalize the definition to
the twisted case.

Following Segal and Quillen's ideas in algebraic K-theory we can
do better in the untwisted case. Consider the category
$\widehat{\CC}$ of virtual objects of $\CC$, namely the objects of
$\widehat{\CC}$ are pairs of vector spaces $(V_0,V_1)$ and so
that a morphism from $(V_0,V_1)$ to $(U_0,U_1)$ is an equivalence
of a triple $[W;f_0,f_1]$ where $W$ is a vector space and
$f_i\colon V_i \oplus W \to U_i$ is an isomorphism. We say that
$(W;f_0,f_1)\sim(W',f_0',f_1')$ if and only if there is an
isomorphism $g\colon W \to W'$ such that $f_i' \circ (
{\mathrm{id}}_{V_i} \oplus g ) = f_i$.

It is a theorem of Segal that $B\widehat{\CC}$ is homotopy
equivalent to the space of Fredholm operators $\FF(\HH)$. But
while for a finite group $G$ it would be wrong to define $K(G)$
as $[BG,\FF(\HH)]$ it is still correct to say that $K(G)$ is the set
of isomorphism classes of functors $\overline{G} \to
\widehat{\CC}$. We can similarly define the $K$-theory of an
orbifold given by a groupoid $\GG$ by functors of the form $\GG
\to \widehat{\CC}$.

Let us again consider the case of a smooth manifold $M$. With
this in mind we would like to have a group model for the space
$\FF$ of Fredholm operators. One possible candidate is the
following.

\begin{definition}\cite{PressleySegal} For a given Hilbert space $\HH$
by a \emph{polarization of $\HH$} we mean a decomposition
$$\HH =\HH_+ \oplus \HH_- $$ where $\HH_+$ is a complete infinite
dimensional subspace of $\HH$ and $\HH_-$ is its orthogonal
complement.

We define the group $\GLres$ to be the subgroup of
${\mathrm{GL}}(\HH)$ consisting of operators $A$ that when written
with respect to the polarization $\HH_+ \oplus \HH_-$ look like
$$A = \left(
      \begin{array}{cc}
      a & b \\
      c & d
      \end{array}
      \right)
$$
where $a\colon \HH_+ \to \HH_+$ and $d \colon \HH_- \to \HH_-$
are Fredholm operators, and $b \colon \HH_- \to \HH_+$ and $c
\colon \HH_+ \to \HH_-$ are Hilbert-Schmidt operators.
\end{definition}

We have the following fact.

\begin{proposition} The map $ \GLres \to \FF \colon A \mapsto a $
is a homotopy equivalence. Therefore $K(M)=[M,\GLres]$.
\end{proposition}

Consider a gerbe $\LL$ with characteristic class $\alpha$ as a
map $ M \to BB\U{1} = B \Proj \U{\HH} $, then we get a Hilbert
projective bundle $Z_\alpha(M) \to M$. Then we form a
$\GLres$-principal bundle over $M$ as follows. We know
\cite{CohenJonesSegal} that polarized Hilbert bundles over $M$
are classified by its characteristic class in $K^1(M)$, for in
view of the Bott periodicity theorem such bundles are classified
by maps $M \to B\GLres = B  B U  = U $, namely by elements in
$K^1(M)$. This produces the desired map $\Gerbe(M) = [M , BBU(1)]
\to [M, U] = K^1(M)$. In several applications it is easier to
start detecting gerbes by means of their image under this map (in the
smooth case, the relation to gerbes and quantum field theory of the  
 $\GLres$-bundles can be found in \cite{Mickelsson}).

\subsection{The twisted theory.}

In this section we are going to ``twist'' vector bundles via gerbes. So for $\RR
\twoarrows \UU$ groupoid associated to the orbifold $X$ and $\LL$ a gerbe over $\RR$

\begin{definition}\label{ourtwist}
An $n$-dimensional $\LL$-twisted bundle over $\RR \twoarrows \UU$ is a groupoid
extension  of it, $\RR \times \complex^n \twoarrows \UU \times \complex^n$ and a
function $\rho : \RR \to \GLC{n}$ such that
$$i^* \rho = \rho^{-1} \ \ \ \ \ \& \ \ \ \ \ (\pi_1^* \rho) \circ (\pi_2^* \rho) \circ ((im)^* \rho)= \theta_\LL \cdot Id_{\GLC{n}}$$
where  $\theta_\LL : \RR \times \RR \to \U{1}$ is the trivialization of triple
intersection $( \pi_1^* \LL \cdot \pi_2^* \LL \cdot (im)^* \LL)
\stackrel{\theta_\LL}{\cong} 1$,
 $Id_{\GLC{n}}$ is the identity of $\GLC{n}$ and the functions $\tilde{s}, \tilde{t}, \tilde{e}, \tilde{i}, \tilde{m}$
are defined in the same way as for bundles.
\end{definition}

We have the following equivalent definition.

\begin{proposition}\label{twiinvar}
   To have an $n$-dimensional $\LL$-twisted bundle over $\RR \twoarrows
   \UU$ is the same thing as to have a vector bundle $E \to \UU$
   together with a given isomorphism $$\LL \otimes t^* E \cong s^* E $$
\end{proposition}

Notice that we then have a canonical isomorphism
$$ m^* \LL \otimes \pi_2^*t^* E \cong \pi_1^* \LL \otimes \pi_2 ^* \LL
\otimes \pi_2^* t^* E \cong \pi_1^* \LL \otimes \pi_2^*(\LL \otimes
t^* E) \cong \pi_1 ^* \LL \otimes \pi_2^* s^* E$$

We can define the corresponding Whitney sum of $\LL$-twisted bundles, so for and
$n$-dimensional $\LL$-twisted bundle with function $\rho_1 :\RR \to \GLC{n}$ and for
an $m$-dimensional one with $\rho_2 :\RR \to \GLC{m}$, we can define a groupoid
extension $\RR \times \complex^{n+m} \twoarrows \UU \times \complex^{n+m}$ with
function $\rho :\RR \to \GLC{n+m}$ such that:
\begin{displaymath}
\rho(g) = \left(\begin{array}{c|c}
\rho_1(g) & 0\\
\hline
0 & \rho_2(g)
\end{array} \right) \in \GLC{n+m}
\end{displaymath}

\begin{definition}
The Grothendieck group generated by the isomorphism classes of $\LL$ twisted bundles
over the orbifold $X$ together with the addition operation just defined is called the
$\LL$ twisted $K$-theory of $X$ and is denoted by $\Kgpd{\LL}(X)$.
\end{definition}

Moerdijk and Pronk \cite{MoerdijkPronk} proved that the isomorphism classes of
orbifolds are in 1-1 correspondence with the classes of \'etale, proper groupoids up
to Morita equivalence. The following is a direct consequence of the definitions.
\begin{lemma}
The construction of $\Kgpd{\LL}(X)$ is independent of the groupoid that is associated
to $X$
\end{lemma}

Similarly as we did with bundles  over groupoids in lemma \ref{isomorphicbundles}, we
can determine when two $\LL$-twisted bundles are isomorphic
\begin{proposition}
An isomorphism of $\LL$-twisted bundles  over $\RR \twoarrows \UU $ (with maps $\rho_i : \RR \to \GLC{n}$ for $i=1,2$)
is determined by a map $\delta : \RR \to \GLC{n}$
such that
\begin{displaymath}
\begin{array}{ccc}
\RR \times \complex^n & \stackrel{\Psi}{\To} &\RR \times \complex^n\\
(r, \xi) & \mapsto  & (r, \delta(r)\xi)
\end{array}
\ \ \ \ \ \
\begin{array}{ccc}
\UU \times \complex^n & \stackrel{\psi}{\To} &\UU \times \complex^n\\
(u, \xi) & \mapsto  & (u, \delta(e(u))\xi)
\end{array}
\end{displaymath}
satisfying $\delta(i(r)) \rho_1(r) = \rho_2(r) \delta (r)$ and $\delta(r) = \delta (e  s(r))$.
\end{proposition}
\begin{proof}
The proof is the same as in lemma \ref{isomorphicbundles}. For $(r,r') \in \RR \timests \RR$ we get:
\begin{eqnarray*}
\theta_\LL id_{\GLC{n}} & = & \rho_1(im(r,r')) \rho_1(r') \rho_1(r) \\
& = & \Big( \delta(m(r,r'))^{-1} \rho_2(im(r,r')) \delta(im(r,r')) \Big)\\
& & \ \ \ \ \ \ \ \ \ \ \ \ \ \ \ \ \ \ \ \ \ \
 \Big( \delta(i(r'))^{-1} \rho_2(r') \delta(r') \Big) \Big( \delta(i(r))^{-1} \rho_2(r) \delta(r) \Big) \\
& = &  \delta(m(r,r'))^{-1} \Big( \rho_2(m(r,r')) \rho_2(r') \rho_2(r) \Big) \delta(r) \\
& = & \delta(m(r,r'))^{-1} \Big( \theta_\LL id_{\GLC{n}} \Big) \delta(r) \\
& = & \theta_\LL id_{\GLC{n}}
\end{eqnarray*}
\end{proof}
Using the group structure of $\Gerbe(\RR \twoarrows \UU)$
we can define a product between bundles twisted by different gerbes, so for
$\LL_1$ and $\LL_2$ gerbes over $X$
$$\Kgpd{\LL_1}(X) \otimes  \Kgpd{\LL_2}(X) \to \Kgpd{\LL_1 \otimes \LL_2}(X)$$
$$(\RR \times \complex^n, \rho_1) \otimes (\RR \times \complex^m, \rho_2) \mapsto (\RR \times \complex^{mn}, \rho_1 \otimes \rho_2))$$
which is well defined because
\begin{eqnarray*}
\lefteqn{(i  m)^*( \rho_1 \otimes \rho_2) \circ \pi_2^* (\rho_1 \otimes \rho_2) \circ \pi_1^*( \rho_1 \otimes \rho_2)} & & \\
\ \ \ & = & \left( (im)^* (\LL_1 \otimes \LL_2) \cdot \pi_2^* (\LL_1 \otimes \LL_2) \cdot \pi_1^* (\LL_1 \otimes \LL_2 \right) \cdot Id_{\GLC{mn}}\\
\ \ \  & = & \theta_{\LL_1 \otimes \LL_2} Id_{\GLC{nm}}\\
\ \ \  & = & \theta_{\LL_1} Id_{\GLC{n}} \otimes \theta_{\LL_2} Id_{\GLC{m}}\\
\ \ \ & = & \left( ((im)^* \LL_1 \cdot \pi_2^* \LL_1 \cdot \pi_1^* \LL_1 ) \cdot Id_{\GLC{n}} \right) \otimes \left( ((im)^* \LL_2 \cdot \pi_2^* \LL_2 \cdot \pi_1^* \LL_2) \cdot Id_{\GLC{m}} \right) \\
\ \ \ & = & \left( ((i  m)^* \rho_1) \circ (\pi_2^* \rho_1) \circ (\pi_1^* \rho_1)
\right) \otimes \left( ((i  m)^* \rho_2) \circ (\pi_2^* \rho_2) \circ (\pi_1^*
\rho_2) \right)
\end{eqnarray*}

and we can define the total twisted orbifold $K$-theory of $X$ as
$$T \Kgpd{}(X) = \bigoplus_{\LL \in \Gerbe(\RR \twoarrows \UU)} \Kgpd{\LL}(X)$$
This has a ring structure due to the following proposition.

\begin{proposition}\label{muygeneral}
   The twisted groups $\Kgpd{\LL}(\GG)$ satisfy the following properties:
   \begin{itemize}
   \item[1.] If $\langle \LL \rangle =0$ then $\Kgpd{\LL}(\GG) = \Kgpd{}(\GG)$ in
   particular if $\GG$ represent the orbifold $X$ then $\Kgpd{\LL}(\GG) = \Korb{}(X)$.
   \item[2.] $\Kgpd{\LL}(\GG)$ is a module over $\Kgpd{}(\GG)$
   \item[3.] If $\LL_1$ and $\LL_2$ are two gerbes over $\GG$ then there is a
   homomorphism
   $$ \Kgpd{\LL_2}(\GG) \otimes \Kgpd{\LL_2}(\GG) \longrightarrow
   \Kgpd{\LL_1\otimes\LL_2}(\GG)$$
   \item[4.] If $\psi\colon \GG_1 \longrightarrow \GG_2$ is a groupoid homomorphism
   then there is an induced homomorphism
   $$ \Kgpd{\LL}(\GG_2) \longrightarrow \Kgpd{\psi^* \LL}(\GG_1) $$
   \end{itemize}
\end{proposition}

\begin{example}
In  the case when $Y$ is the orbifold universal cover of $X$ with orbifold
fundamental group $\pi_1^{orb}(X) =H$,
 we  can take a discrete torsion $\alpha \in H^2(H,\U{1})$ and define the twisted $K$-theory
of $X$ as in Definition \ref{alfatwistedk}. Let's associate to $X$ the groupoid
constructed in example \ref{exampleorbifolduniversalcover}; we want to construct a
gerbe $\LL$  over $\RR_Y \times H \twoarrows \UU_Y$ so that the twisted
$\Kgpd{\LL}(X)$ is the same as the twisted $\Korb{\alpha}^{AR}(X)$ of section
\ref{subsectionorbifoldKtheory} (we added the upperscipts ${}^{AR}$ to denote that is
the twisted $K$ theory defined by A. Adem and Y. Ruan).

The discrete torsion $\alpha$ defines a central extension of $H$
$$ 1 \To \U{1} \To \widetilde{H} \To H \To 1$$
and doing the cartesian product with $\RR_Y$ we get a line bundle

\begin{displaymath}
\begin{array}{ccccc}
\U{1} & & \To &\LL_\alpha = & \RR_Y \times \widetilde{H} \\
 & & & & \downarrow \\
 & & & & \RR_Y \times H
\end{array}
\end{displaymath}
which, by lemma \ref{gerbesdiscretegroup} and the fact that the
line bundle structure comes from the lifting of $H$, becomes a
gerbe over $\RR_Y \times H \twoarrows \UU_Y$ . Clearly this gerbe
only depends on the one defined in lemma
\ref{gerbesdiscretegroup} for the group $H$. For $E \to X$ an
$\alpha$-twisted bundle over $Y$, an element of
$\Korb{\alpha}^{AR}(X)$, comes with an action of $\widetilde{H}$
in $E$ such that it lifts the one of $H$ in $Y$; choosing specific
lifts $\tilde{h},\tilde{g} \in \widetilde{H}$ for every  $h,g \in
H$, and $e \in E$, we have $\tilde{g} (\tilde{h} (e))=
\alpha(g,h) \widetilde{gh}(e)$. As $E$ is a bundle over $Y$, it
defined by a map $\rho:\RR_Y \to \GLC{n}$, and for $h \in H$, it
defines an isomorphism $\RR_Y \times \complex^n \stackrel{h}{\to}
\RR_Y \times \complex^n$, with $\eta_h : \RR_Y \to \GLC{n}$
 such that $(r,\xi) \stackrel{h}{\mapsto} (hr, \eta_h(r)\xi)$.
The $\LL_\alpha$-twisted bundle over $\RR_Y \times H \twoarrows
\UU_Y$ that $E$ determines, is given by the groupoid $\RR_Y
\times H \times \complex^n \twoarrows \UU_Y \times \complex^n$
and the map
\begin{eqnarray*}
\delta : \RR_Y \times H & \to & \GLC{n} \\
(r,h) & \mapsto & \rho(hr) \eta_h(r)
\end{eqnarray*}
Because $h$ is an isomorphism of groupoids
$\RR_Y\times\complex^n\to\RR_Y\times\complex^n$, it commutes with
the source and target maps, and this in turns implies that for
$h\in H$ and $r \in \RR_Y$ we have
$\eta_h(r)=\eta_h(e_Y(s_Y(r)))$ and
$\eta_h(r)\rho(r)=\rho(hr)\eta_h(r)$. In order to prove that this
bundle is $\LL_\alpha$-twisted it is enough to check that the
multiplication satisfies the specified conditions.

We will make use of the following diagram in the calculation
        $$\xymatrix{
        x \ar[r]^r \ar[dr]_w & y \ar[d]^v\\
         & z}
        \ \ \ \  \stackrel{h}{\longrightarrow} \ \ \ \
        \xymatrix{
        x' \ar[r]^{r'} \ar[dr]_{w'} & y' \ar[d]^{v'}\\
         & z'}
        \ \ \ \ \stackrel{j}{\longrightarrow} \ \ \ \
        \xymatrix{
        x'' \ar[r]^{r''} \ar[dr]_{w''} & y'' \ar[d]^{v''}\\
         & z''}
        $$
where the $r$'s, $v$'s and $w$'s belong to $\RR_Y$, the $x$'s, $y$'s and $z$'s belong
to $\UU_Y$ and $h,j\in H$. We have that $m_Y(r,v)=w$ and
$$ m((r,h),(v',j))=(m(r,v),jh)=(w,jh).$$

Also
\begin{eqnarray*}
\delta(v',j)\delta(r,h) & = & \rho(j v')\eta_j(v')\rho(hr)\eta_h(r)\\
                          & = & \eta_j(v')\rho(v')\rho(r')\eta_h(r) \\
                          & = & \eta_j(v')\rho(w')\eta_h(r)\\
                          & = & \eta_j(v')\rho(hw)\eta_h(w)\\
                          & = & \eta_j(v')\eta_h(w)\rho(w)\\
                          & = & \eta_j(v')\eta_h(r)\rho(w)\\
                          & = & \alpha(j,h)\eta_{jh}(w)\rho(w)\\
                          & = & \alpha(j,h)\delta(w,j h)
\end{eqnarray*}

The fact that $\eta_j(v')\eta_h(r)=\alpha(j,h)\eta_{jh}(w)$ is precisely the fact
that $E$ is an $\alpha$-twisted bundle. This proves that $(\RR_Y\times H \times
\complex^n \rightrightarrows \UU_Y \times \complex^n, \delta)$ is endowed with the
structure of a $\LL_\alpha$-twisted bundle over $\RR_Y \times H \rightrightarrows
\UU_Y$.

\end{example}

Conversely, if we have the $\LL_\alpha$-twisted bundle, it is clear how to obtain the
maps $\rho$ and $\eta_h$. From the previous construction we can see that when $h =
id_H$, the map $\delta$ determines $\rho$ (i.e. $\delta(r,id_H)=\rho(r)$); hence
$\eta_h(r) = \delta(r,h) \rho(hr)^{-1}$. Thus we can conclude,

\begin{theorem}\label{oursadem} In the above example
$\Kgpd{\LL_\alpha}(X) \cong \Korb{\alpha}^{AR}(X)$
\end{theorem}

From \ref{Poincare}, \ref{smoothtwist} and \ref{character} we have

\begin{proposition} If $X$ is a smooth manifold and the characteristic class of the
gerbe $\LL$ is the torsion element $[H]$ in $H^3(M,\integer)$ then
$$\Kgpd{\LL}(X) \cong K_{[H]}(X)$$
\end{proposition}

It remains to verify that the twisting  $\Kgpd{\LL}(X)$ coincides
with the twisting $K_\alpha (X)$ defined by \ref{AtiyahTwist}.

\begin{proposition} Whenever $\alpha = \langle \LL
\rangle $ is a torsion class and $X=M$ is a smooth manifold then
$\Kgpd{\LL}(M)= K_\alpha (M)$.
\end{proposition}
\begin{proof} We will use the following facts.
   \begin{theorem}[Serre \cite{DonovanKaroubi}] Let $M$ be a CW-complex.
   If a class $\alpha \in H^3(M,\integer)$ is a torsion element
   then there exists a principal bundle $Z\to M$ with structure
   group $\Proj\U{n}$ so that when seen as an element
   $ \beta \in [M,B\Proj\U{n}] \to [M,B\Proj{\mathrm{U}}] = [M , B B \U{1}] = [M , B
   K(\integer,2)] = [M, K(\integer, 3)]= H^3(M,\integer)$ then
   $\alpha = \beta$. In other words, the image of $[M,B\Proj\U{n}]
   \to H^3(M,\integer)$ is exactly the subgroup of torsion elements.
   \end{theorem}
   \begin{theorem}[Segal \cite{SegalKFred}]\label{SegalKF} Let $\HH$ be a $G$-Hilbert space
   in which every irreducible representation of $G$ appears infinitely many times.
   Then the equivariant
   index map $$ {\mathrm{ind}}_G \colon [Z , \FF]_G  \longrightarrow K_G(Z) ,$$
   is an isomorphism (where $G$ acts on $\FF$ by conjugation.)
   \end{theorem}
   \begin{lemma}
   Let $X=Z/G$ be an orbifold where
   the Lie group $G$ acts on $Z$. Let $\alpha \in H^2(G,\U{1})$ define a group extension
   $1 \to \U{1} \to \tilde{G} \to G \to 1$. Consider the natural
   homomorphism $ \psi \colon K_{\tilde{G}}(Z) \to K_{\tilde{G}}(G) \cong
   R(\U{1})$ - it can also be seen as the composition of
  $  K_{\tilde{G}}(Z) \to K_{\tilde{G}}(*) \cong R(\tilde{G}) \to R(\U{1})$.
   Let $\Delta : \U{1} \to \U{n}$ be the diagonal
   embedding representation. Then $$\Korb{\alpha}^{AR}(X) = \psi
   ^{-1} (\Delta).$$
   \begin{proof} The orbifold $X$ is represented by the groupoid
   $ \GG =( Z \times G \rightrightarrows Z )$, while the gerbe $\LL_\alpha$
   is represented by the central extension of groupoids ( proposition \ref{gerbeextension} )
   $$  1 \to \overline{\U{1}} \to \widetilde{\GG} \to \GG \to 1$$
   where $\widetilde{\GG} = ( Z \times \tilde{G}  \rightrightarrows Z)$.
   Therefore, using the fact that $\Kgpd{}(\overline{\U{1}}) =
   R(U(1))$ we get the surjective map
   $$ \Kgpd{}(\widetilde{\GG}) \to R(\U{1}), $$
   using \ref{ourtwist} and observing that $ \Kgpd{}(\widetilde{\GG})
   = K_{\tilde{G}}(Z)$ we get the result.
   \end{proof}
   \end{lemma}
Let us consider in the previous lemma the situation where $X=M$
is smooth, $Z$ is Serre's principal $\Proj\U{n}$-bundle
associated to $\alpha=\langle\LL\rangle$,  $G=\Proj\U{n}$ and
$\beta$ is the class in $H^2(\Proj\U{n},\U{1})$ labeling the
extension
$$ 1 \to \U{1} \to \U{n} \to \Proj\U{n} \to 1 .$$ Then using
Theorem \ref{oursadem} and \ref{SegalKF} we get that
$$ \Kgpd{\LL}(Z/G) = \Korb{\beta}^{AR}(Z/G) = \psi^{-1}(\Delta) \subseteq
K_{\U{n}} (Z) = [Z, \FF]_{\U{n}}.$$

Notice that by \ref{AtiyahTwist} $K_\alpha(M)$ is defined as the
homotopy classes of sections of the bundle $\FF_\alpha = Z
\times_{\Proj\U{n}} \FF$. This space of sections can readily be
identified with the space $[Z , \FF_\alpha]_{\Proj\U{n}}$ and the
proposition follows from this.
\end{proof}

We should point out here that the theory so far described is
essentially empty whenever the characteristic class $\langle \LL
\rangle$ is a non-torsion element in $H^3(M,\integer)$. The
following is true.

\begin{proposition} If there is an $n$-dimensional $\LL$-twisted bundle
over the groupoid $\GG$ then $\langle\LL\rangle ^n = 1 $.
\end{proposition}
\begin{proof}
Consider the equations,
$$i^* \rho = \rho^{-1} \ \ \ \ \ \& \ \ \ \ \ (\pi_1^* \rho) \circ (\pi_2^* \rho) \circ ((im)^* \rho)= \theta_\LL \cdot Id_{\GLC{n}}$$
and take determinants in both equations, we get
$$i^* \det \rho = \det \rho^{-1} \ \ \ \ \ \& \ \ \ \ \ (\pi_1^* \det \rho) \circ (\pi_2^* \det \rho) \circ ((im)^* \det \rho)= \det \theta_\LL \cdot Id_{\GLC{n}}$$
defining $f = \det \rho$ we have
$$i^* f = f^{-1} \ \ \ \ \ \& \ \ \ \ \ (\pi_1^* f) \circ (\pi_2^* f) \circ ((im)^* f)= \theta_\LL ^ n$$
this means that the coboundary of $f$ is $\theta^n$. This
concludes the proof.
\end{proof}

Another way to think of this is by noticing that  if we restrict
the central extension
$$ 1\to \U{1} \to \U{n} \to \Proj\U{n} \to 1$$
to the subgroup $S\U{n}$ we get the $n$-fold covering map
$$ 1\to \integer_n \to S\U{n} \to \Proj\U{n} \to 1 $$
where the kernel $\integer_n$ is the group of $n$-roots of unity.

In any case we need to consider a more general definition when
the class $\langle \LL \rangle$ is a non-torsion class.

An obvious generalization of \ref{AtiyahTwist} would be to
consider the class of the gerbe $\alpha = \langle\LL\rangle \in
H^3(B\GG,\integer)$ and consider $K_\alpha(B\GG)$ in the sense of
\ref{AtiyahTwist}. This works well for a manifold, but
unfortunately for a finite group and the trivial gerbe $\alpha=1$
we have that $K_\alpha(B \overline{G}) = R(G)\widehat{\ }$ and
not $R(G)$ as we should have (this is exactly the problem we
encountered in the last section with $[BG,\FF(\HH)]$).

Fortunately one of the several equivalent definitions of
\cite{Murray2} can be carefully generalized to serve our purposes.
To motivate this definition consider the following situation.

Suppose first that the class $\alpha = \langle\LL\rangle \in
H^3(B\GG,\integer)$ is a torsion class. Take any $\alpha$-twisted
vector bundle $\rho$ so that
$$i^* \rho = \rho^{-1} \ \ \ \ \ \& \ \ \ \ \ 
(\pi_1^* \rho) \circ (\pi_2^* \rho) \circ ((im)^* \rho)= \theta_\LL \cdot Id_{\U{n}}$$
and let $\beta \colon \GG \to \Proj\U{n}$  be the projectivization
of $\rho \colon \GG \to \U{n}$. Then $\beta$ is a \emph{bona fide}
groupoid homomorphism, in other words
$$i^* \beta = \beta^{-1} \ \ \ \ \ \& \ \ \ \ \
 (\pi_1^* \beta) \circ (\pi_2^* \beta) \circ ((im)^* \beta)= Id_{\Proj\U{n}}$$
as equations in $\Proj\U{n}$. Then $\alpha$ as a map $B\GG \to
B\Proj \U{\HH}$ is simply obtained as the realization of the composition of  $\beta \colon \GG
\to \Proj\U{n}$  with the natural inclusion $\Proj\U{n}
\hookrightarrow \Proj\U{\HH}$. Fix now and for all a $\rho_0$ and
a $\beta_0$ constructed in this way.

Define the semidirect product $\Semiproductofin$ as the group
whose elements are the pairs $(S,T)$ where $S\in \U{n}$ and
$T\in\Proj\U{n}$ with multiplication
$$(S_1,T_1)\cdot(S_2,T_2) = (S_1 T_1 S_2 T_1^{-1} , T_1 T_2).$$

Consider the family of groupoid homomorphisms $f\colon \GG \to
\Semiproductofin$ that make the following diagram commutative

$$
        \xymatrix{
        \GG \ar[r]^f \ar[rd]_{\beta_0} & \Semiproductofin \ar[d]^{q_2}\\
          & \Proj\U{n}}
$$
where $q_1$ is the projection onto $\U{n}$ and $q_2$ the
projection onto $\Proj \U{n}$.

Given a homomorphism $f \colon \GG \to \Semiproductofin$ like
above we can then write $\rho = (q_1 \circ f) \cdot \rho_0$ and
verify that $\rho$ satisfies the conditions to define a twisted
vector bundle over $\GG$. Conversely given a twisted vector
bundle $\rho$ we can define a homomorphism $f$ by means of the
formula $$ f(g) = ( \rho(g) \rho_0(g)^{-1} , \beta_0(g)).$$
Therefore in the case of a torsion class $\alpha$ a homomorphism
$f \colon \GG \to \Semiproductofin $ so that $q_1 f = \beta_0 $
is another way of encoding a twisted vector bundle.

In the case of a non-torsion class $\alpha$ we need to consider
infinite dimensional vector spaces. So we let $\KK$ be the space
of compact operators of a Hilbert space $\HH$. Let us write $\UK$
to denote the subgroup of $\U{\HH}$ consisting of unitary
operators of the form $I+K$ where $I$ is the identity operator
and $K$ is in $\KK$. If $h \in \Proj\U{\HH}$ and $g \in \UK$ then
$hgh^{-1} \in \UK$ and therefore we can define $\Semiproducto$.
We can define now the $K$-theory for an orbifold $X$ given by
$\GG$ twisted by a gerbe $\LL$ with non-torsion class $\alpha
\colon \GG \to \Proj \U({\HH})$ (cf. \ref{GerbeAsHomo}).

\begin{definition} The set of isomorphism classes of groupoid
homomorphisms $f\colon \GG \to \Semiproducto$ that make the
following diagram commutative
$$
        \xymatrix{
        \GG \ar[r]^f \ar[rd]_{\alpha} & \Semiproducto \ar[d]^{q_2}\\
          & \Proj\U{\HH}}
$$
is $\Kgpd{\LL}(\GG)$ the groupoid $K$-theory of $\GG$ twisted by
$\LL$.
\end{definition}

This definition works for a gerbe whose class is non-torsion and
has the obvious naturality conditions. In particular it becomes
\ref{AtiyahTwist} if the groupoid represents a smooth manifold.
The discussion immediately before the definition shows that this
definition generalizes the one given before for $
\Kgpd{\LL}{\GG}$ when $\langle\LL\rangle$ was torsion. Proposition
\ref{muygeneral} remains valid.

In view of theorem \ref{holonomy} and theorem \ref{oursadem} we
can reformulate theorem \ref{ADEMRUANTH} as follows.

\begin{theorem}
Let $\Xx$ be a Leray groupoid representing an orbifold $X/\Gamma$
with $\Gamma$ finite , $\LL$ a gerbe over $\Xx$ coming from
discrete torsion, and let $\Lambda$ be the holonomy inner local
system defined in \ref{holonomy}. Then
$$\Kgpd{\LL}(\Xx) \otimes \complex \cong H^*_{orb,\Lambda}(\Xx; \complex)$$
\end{theorem}

It is natural to conjecture that the previous theorem remains
true even in the gerbe $\LL$ is arbitrary and $\Xx$ is any proper
\'etale groupoid. We will revisit this issue elsewhere.

The following astonishing result of Freed, Hopkins and Teleman
can be written in terms of the twisting described in this
section. For more on this see \cite{LupercioUribe2}.

\begin{example} \cite{Freed}
Let $G$ be connected, simply connected and simple. Consider the groupoid $\GG=(G
\times G \rightrightarrows G)$ where $G$ is acting on $G$ by \emph{conjugation}. Let
$h$ be the dual Coxeter number of $G$. Let $\LL$ be the gerbe over $\GG$ with
characteristic class $ \dim(G) + k + h \in H^3_G(G) $. Then
$$ \Kgpd{\LL}(\GG) \cong V_k(G)$$
where $V_k(G)$ is the Verlinde algebra at level $k$ of $G$.
\end{example}

\subsection{Murray's bundle gerbes.}

The theory described in the previous sections is interesting even
in the case in which the orbifold $X$ is actually a smooth
manifold $M=X$. In this case M. Murray et. al. \cite{Murray1,
Murray2} have recently proposed a way to interpret the twisted
$K$-theory ${}^L K(M)$ in terms of bundle gerbes. Bundle gerbes
are geometric objects constructed on $M$ that give a concrete
model for a gerbe over $M$ \cite{Murray0}. The purpose of this
section is to explain how the theory of bundle gerbes can be
understood in terms of groupoids.

\begin{definition}
   A bundle gerbe over $M$ is a pair $(L,Y)$ where $Y \stackrel{\pi}{\longrightarrow}
   M$ is a surjective submersion and $L\stackrel{p}{\longrightarrow}Y\: {}_{\pi}  \! \times_{\pi} Y =
   Y^{[2]}$ is a line bundle satisfying
   \begin{itemize}
    \item $L_{(y,y)} \cong \complex$
    \item $L_{(y_1,y_2)} \cong L_{(y_2,y_1)}^*$
    \item $L_{(y_1,y_2)} \otimes L_{(y_2,y_3)} \cong L_{(y_1,y_3)}$
   \end{itemize}
\end{definition}

We start the translation to the groupoid language with the following definition.

\begin{definition}
   Given a manifold $M$ and a surjective submersion $Y \stackrel{\pi}{\longrightarrow}
   M$ we define the groupoid $\GG(Y,M) = ( \RR \rightrightarrows \UU ) $ by
   \begin{itemize}
      \item $\RR=Y^{[2]} = Y\: {}_{\pi}  \! \times_{\pi} Y $
      \item $\UU = Y$
      \item $s = p_1 : Y^{[2]} \to Y$, $s(y_1,y_2)=y_1$ and $t= p_2 : Y^{[2]} \to Y$,
      $t(y_1,y_2)=y_2$
      \item $m((y_1,y_2),(y_2,y_3))=(y_1,y_3)$
   \end{itemize}
\end{definition}

From Definition \ref{gerbedef} we immediately obtain the following.

\begin{proposition}
   A bundle gerbe (L,Y) over $M$ is the same as a gerbe over the groupoid $\GG(Y,M)$.
\end{proposition}

We will write $\LL(L,Y)$ to denote the gerbe over $\GG(Y,M)$ associated to the bundle
gerbe $(L,Y)$. Notice that the groupoid $\GG(Y,M)$ is not necessarily \'etale, but it
is Morita equivalent to an \'etale groupoid. Let $\MM(M,U_\alpha)$ be the \'etale
groupoid associated to a cover $\{U_\alpha\}$ of $M$ as in \ref{examplemanifold}.

\begin{proposition}
   The groupoid $\GG(Y,M)$ is Morita equivalent to $\MM(M,U_\alpha)$ for any open cover
   $\{U_\alpha\}$ of $M$.
\end{proposition}
\begin{proof} Since all groupoids $\MM(M,U_\alpha)$ are Morita equivalent (for any two
open covers have a common refinement,) it is enough to consider the groupoid $\MM(M,M)
= (M \rightrightarrows M) $ coming from the cover consisting of one open set. The
source and target maps of $\MM(M,M)$ are both identity maps. Then the proposition
follows from the fact that the following diagram is a fibered square
        $$\xymatrix{
        Y^{[2]} \ar[r]^{\pi s} \ar[d]_{(s,t)} & M \ar[d]^{s \times t}\\
        Y \times Y \ar[r]^{\pi \times \pi} & M \times M}
        $$

\end{proof}

\begin{cor}\label{weisomurray} The group of bundle gerbes over $M$ is isomorphic to the
group $\Gerbe(\MM(M,U_\alpha))$ for the Leray groupoid $\MM(M,U_\alpha)$ representing
$M$. In particular there is a bundle gerbe in every Morita equivalence class of
gerbes over $M$.
\end{cor}

Murray \cite{Murray1} defines a characteristic class for a bundle gerbe $(L,Y)$ over
$M$ as follows.

\begin{definition}
    The \emph{Dixmier-Douady class} of $d(P)=d(P,Y)\in H^3(M,\integer)$ is defined as
    follows. Choose a Leray open cover $\{ U_\alpha\}$ of $M$. Choose sections
    $s_\alpha \colon U_\alpha \to  Y$, inducing  $(s_\alpha,s_\beta)\colon
    U_\alpha \cap  U_\beta \to Y^{[2]}$. Choose sections $\sigma_{\alpha \beta}$ of
    $(s_\alpha,s_\beta)^{-1}(P)$ over $ U_\alpha \cap  U_\beta  $. Define
    $g_{\alpha \beta \gamma} \colon U_\alpha \cap U_\beta \cap U_\gamma \to
    \complex^\times$ by $$ \sigma_{\alpha\beta} \sigma_{\beta\gamma} =
    \sigma_{\alpha\gamma} g_{\alpha \beta \gamma} $$
    Then $d(L)=[g_{\alpha \beta \gamma}] \in H^2(M,\complex^\times) \cong
    H^3(M,\integer)$.
\end{definition}

Then in view of propositions \ref{character} and \ref{indepchar} we have the following.

\begin{proposition}\label{ddmurray}
 The Dixmier-Douady class $d(L,Y)$ is equal to the characteristic class $\langle
 \LL(L,Y) \rangle$ defined above \ref{character}. Moreover the assignment $(L,Y)
 \mapsto {g_{\alpha \beta \gamma}} $ realizes the isomorphism of \ref{weisomurray}.
\end{proposition}

\begin{definition}\cite{Murray1} A bundle gerbe $(L,Y)$ is said to be \emph{trivial}
whenever $d(L,Y)=0$. Two bundle gerbes $(P,Y)$ and $(Q,Z)$ are called \emph{stably
isomorphic} if there are trivial bundle gerbes $T_1$ and $T_2$ such that
$$ P\otimes T_1 \simeq Q\otimes T_2. $$
\end{definition}

The following is an easy consequence of \ref{character}, \ref{indepchar} and
\ref{ddmurray}.

\begin{lemma} The Dixmier-Douady class is a homomorphism from the group of bundle
gerbes over $M$ with the operation of tensor product, and $H^3(M,\integer)$
\end{lemma}

\begin{cor} Two bundle gerbes $(P,Y)$ and $(Q,Z)$ are
stably isomorphic if and only if $d(P)=d(Q)$
\end{cor}
\begin{proof} Suppose that $(P,Y)$ and $(Q,Z)$ are
stably isomorphic. Then $ P\otimes T_1 \simeq Q\otimes T_2; $ hence $d( P\otimes T_1 )
= d( Q\otimes T_2)$. Therefore from the previous lemma we have
$d(P)+d(T_1)=d(Q)+d(T_2)$ and by definition of trivial we get $d(P)=d(Q)$.

Conversely if $d(P)=d(Q)$ then $d(P\otimes Q^*)=0$ and then by definition
$T_2=P\otimes Q^*$ is trivial. Define the trivial bundle gerbe $T_1=Q^*\otimes Q$.
Then $ P\otimes T_1 \simeq Q\otimes T_2 $ completing the proof.
\end{proof}

 Given a bundle gerbe $(L,Y)$ over $M$ we will write $\widetilde{\GG}(L,Y,M)$ to
denote the $\overline{\U{1}}$ central groupoid extension of $\GG(Y,M)$ defined by the
associated gerbe, where $$ 1 \to \overline{\U{1}} \to \widetilde{\GG}(L,Y,M) \to
\GG(Y,M) \to 1 .$$ As we have explained before such extensions are classified by
their class in the cohomology group $H^3(B\GG(Y,M),\integer)=H^3(M,\integer)$. As a
consequence of this and \ref{weisomurray} we have.

\begin{theorem}\label{equmurus} Two bundle gerbes $(P,Y)$ and $(Q,Z)$ are
stably isomorphic if and only if $\widetilde{\GG}(P,Y,M)$  is Morita equivalent to
$\widetilde{\GG}(Q,Z,M)$. Therefore there is a one-to-one correspondence between
stably isomorphism classes of bundle gerbes over M and classes in $H^3(M,\integer)$.
The category of bundle gerbes over $M$ with stable isomorphisms is equivalent to the
category of gerbes over $M$ with Morita equivalences.
\end{theorem}

\begin{definition} Let $(L,Y)$ be a bundle gerbe over $M$. We call $(E,L,Y,M)$ a
\emph{bundle gerbe module} if
\begin{itemize}
\item $E\to Y$ is a hermitian vector bundle  over
$Y$
\item We are given an isomorphism $\phi \colon L \otimes
\pi_1^{-1} E \stackrel{\sim}{\longrightarrow} \pi_2^{-1} E $.
\item The compositions $L_{(y_1,y_2)}\otimes(L_{(y_2,y_3)}\otimes E_{y_3})\to
L_{(y_1,y_2)}\otimes E_{y_2} \to E_{y_1}$ and $(L_{(y_1,y_2)}\otimes L_{(y_2,y_3)})
\otimes E_{y_3} \to L_{(y_1,y_3)}\otimes E_{y_3} \to E_{y_1}$ coincide.
\end{itemize}
In this case we also say that the bundle gerbe $(L,Y)$ acts on
$E$. The \emph{bundle gerbe K-theory} $K_{\mathrm{bg}}(M,L)$ is
defined as the Grothendieck group associated to the semigroup of
bundle gerbe modules $(E,L,Y,M)$ for $(L,Y,M)$ fixed.
\end{definition}

As a consequence of \ref{equmurus} and \ref{twiinvar} we have the
following fact.

\begin{theorem} The category of bundle gerbe modules over $(L,Y)$
is equivalent to the category of $\LL(L,Y)$-twisted vector
bundles over $\GG(Y,M)$. Moreover we have
$$ \Kgpd{\LL}(\GG(Y,M)) \cong K_{\mathrm{bg}}(M,L)$$
\end{theorem}

\begin{cor} If the gerbe $\LL$ has a torsion class $[H]$
then $$ K_{[H]}(M) = K_{\mathrm{bg}}(M,L).$$
\end{cor}

\section{Appendix: Stacks, gerbes and groupoids.}

We mentioned at the beginning of Section \ref{Groupoids} that a stack $\XX$ is a space
whose points can carry a ``group valued" multiplicity, and that they are studied by
studying the family $\{\mbox {Hom$(S,\XX)$}\}_S$ where $S$ runs through all possible
spaces (or schemes). In fact by Yoneda's Lemma, as is well known, even when $\XX$ is
an ordinary space that knowing everything for the functor Hom$(\centerdot,\XX)$ is the
same thing as knowing everything about $\XX$. A stack is a category fibered by
groupoids where $\CC_\SSS = \mbox{Hom}(S,\XX)$ with an additional sheaf condition.

A very unfortunate confusion of terminologies occurs here. The word \emph{groupoid}
has two very standard meanings. One has been used along all the previous sections of
this paper. But now we need the second meaning, namely a groupoid is a category where
all morphisms have inverses. In this Appendix we use the word groupoid with both
meanings and we hope that the context is enough to avoid confusion. Both concepts
are, of course, very related.

\subsection{Categories fibered by groupoids}
Let $\CC, \SSS$ be a pair of categories and $p: \SSS \to \CC$ a
functor. For each $U \in Ob(\CC)$ we denote $\SSS_U =p^{-1}(U)$.

\begin{definition} The category $\SSS$ is \emph{ fibered by groupoids} over $\CC$ if

\begin{itemize}
\item For all $\phi : U \to V$ in $\CC$ and $ y \in Ob(\SSS_V)$
there is a morphism $f: x \to y$
in $\SSS$ with $p(f)=\phi$
\item For all $\psi : V \to W$, $\phi: U \to W$, $\chi : U \to V$, $f: x \to y$ and
$g: y \to z$ with $\phi = \psi \circ \chi$, $p(f) = \phi$ and $p(g)=\psi$ there is a
unique $h: x \to z $ such that $f=g \circ h$ and $p(h) = \chi$
\end{itemize}
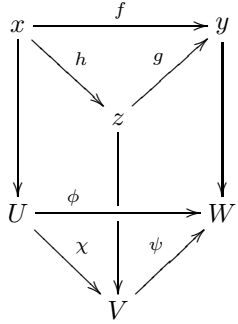
        $$\xymatrix{
        x \ar[rr]^f \ar[dr]^h \ar[dd] & & y \ar[dd]\\
          & z \ar[ur]^g \ar'[d][dd] & \\
        U \ar[rr]^<(.2){\phi} \ar[dr]^\chi & & W\\
         & V \ar[ur]^\psi & }
        $$

\end{definition}

The conditions imply that the existence of the morphism $f:x \to y$ is unique up
 to canonical isomorphism. Then for $\phi :
U \to V$ and $y \in Ob(\SSS_V)$, $f:x \to y$ has been chosen; $x$
will be written as $\phi^*y$ and $\phi^*$ is a functor from
$\SSS_V$ to $\SSS_U$.

\subsection{Sheaves of Categories}

\begin{definition}
A \emph{ Grothendieck Topology} (G.T.) over a category $\CC$ is a
prescription of coverings $ \{ U_\alpha \to U\}_\alpha$ such that:
\begin{itemize}
\item $\{ U_\alpha \to U\}_\alpha$ \& $\{ U_{\alpha \beta} \to U_\alpha\}_\beta$
implies $\{ U_{\alpha \beta} \to U\}_{\alpha \beta}$
\item $\{ U_\alpha \to U\}_\alpha$ \& $V \to U$ implies $\{ U_\alpha \times_U V
 \to V \}_\alpha$
\item $V \stackrel{\cong}{\longrightarrow} U$ isomorphism, implies $\{  V \longrightarrow U\}$
\end{itemize}
A category with a Grothendieck Topology is called a \emph{Site}.
\end{definition}

\begin{example} \label{smallsite} $\CC = Top$, $\{ U_\alpha \to U\}_\alpha$ if $U_\alpha$
is homeomorphic to its image and $U=\bigcup_\alpha im(U_\alpha)$.
\end{example}

\begin{definition}
A \emph{Sheaf} $\FF$ over a site $\CC$ is a functor p:$\FF \to
\CC$ such that
\begin{itemize}
\item For all $ S \in Ob(\CC)$, $ x \in Ob(\FF_S)$ and $f:T \to S \in
Mor(\CC)$ there exists a { unique} $\phi : y \to x \in Mor(\FF)$
such that $p(\phi) = f$
\item For every cover $\{ S_\alpha \to S\}_\alpha$, the following sequence is exact
$$\FF_S \to \prod \FF_{S_\alpha} \twoarrows \prod \FF_{S_\alpha \times_S S_\beta}$$
\end{itemize}
\end{definition}

\begin{definition}
A \emph{Stack in groupoids} over $\CC$ is a functor $p: \SSS \to \CC$ such that
\begin{itemize}
\item $\SSS$ is fibered in groupoids over $\CC$
\item For any $U \in Ob(\CC)$ and $x,y \in Ob(\SSS_U)$, the functor
$${\bf U} \to Sets$$
$$\phi : V \to U \mapsto Hom(\phi^*x, \phi^*y)$$
is a sheaf. $\left(Ob({\bf U})= \{(S,\chi) | S \in Ob(\CC), \chi
\in Hom(S,U) \} \right)$
\item If $\phi_i : V_i \to U$ is a covering family in $\CC$, any descent datum
relative to the $\phi_i$'s, for objects in $\SSS$, is effective.
\end{itemize}
\end{definition}

\begin{example} For $X$ a $G$-set (provided with a $G$ action over it) let $\CC =Top$, the
category of topological spaces, and $S=[X/G]$ the category
defined as follows:
$$Ob([X/G])_S = \{ f: E_S \to X \}$$

the set of all $G$-equivariant maps from principal $G$-bundles
$E_S$ over $S \in Ob(Top)$, and
$$Mor([X/G]) \subseteq Hom_{BG}(E_S,E'_S)$$
given by

       $$\xymatrix{
       E_S \ar@{{<}-{>}}[r] \ar[d]_{(proj,f)} & S \times_{S'} E_{S'} \ar[d]^{1 \times f'}\\
       S \times X \ar@{{<}-{>}}[r] & S \times X}
     $$

With the functor

\begin{eqnarray*}
 p: [X/G] & \to & Top \\
 (f: E_S \to X) & \mapsto & S
\end{eqnarray*}

By definition $[X/G]$ is a category fibered by groupoids, and if the group $G$ is
finite $[X/G]$ is a stack.
\end{example}

\subsection{Gerbes as Stacks}

For simplicity we will start with a smooth $X$ and we will
consider a Grothendieck topology on $X$ induced by the ordinary
topology on $X$ as in \ref{smallsite}. We will follow Brylinski
\cite{Brylinski} very closely. The following definition is
essentially due to Giraud \cite{Giraud}.

\begin{definition} \label{gerbeasstack}
A \emph{gerbe} over $X$ is a sheaf of categories $\CC$ on $X$ so
that
\begin{itemize}
   \item The category $\CC_U$ is a groupoid for every open $U$.

   \item Any two objects $Q$ and $Q'$ of $\CC_U$ are locally
   isomorphic, namely for every $x\in X$ there is a neighborhood
   of $x$ where they are isomorphic.

   \item Every point $x \in X$ has a neighborhood $U$ so that
   $\CC_U$ is non-empty.
\end{itemize}

We will require our gerbes to have as \emph{band} the sheaf
$\Aa=\complex^*$ over $X$. This means that for every open $U
\subset X$ and for every object $Q \in \CC_U$ there is an
isomorphism of sheaves $\alpha \colon
\underline{{\mathrm{Aut}}}(Q) \to \Aa_U$, compatible with
restrictions and commuting with morphisms of $\CC$. Here
$\underline{{\mathrm{Aut}}}(Q)|_V$ is the group of automorphisms
of $P|_V$.
\end{definition}

The relation of this definition to the one we have used is given
by the following

\begin{proposition} \label{gerbesassheaves} To have a gerbe in terms of data
$(\LL_{\alpha\beta})$ as in \ref{gerbessmooth} is the same thing as to have a gerbe
with band $\complex^*$ as a sheaf of categories.
\end{proposition}

\begin{proof}
   Starting with the data in \ref{gerbessmooth} we will construct
   the category $\CC_U$ for a small open set $U$. Since $U$ is
   small we can trivialize the gerbe $\LL|_U$. The objects of
   $\CC_U$ are the set of all possible trivializations
   $(f_{\alpha\beta})$ with the obvious morphisms.

   Conversely suppose that you are given a gerbe as a sheaf of
   categories. Then we construct a cocycle $c_{\alpha\beta\gamma}$
   as in \cite[Prop. 5.2.8]{Brylinski}. We take an object $Q\in
   C_{U_\alpha}$ and an automorphisms $u_{\alpha\beta} \colon
   Q_\alpha |_{U_{\alpha \beta}} \to  Q_\beta |_{U_{\alpha \beta}}$
   and define $h_{\alpha\beta\gamma} = u_{\alpha\gamma}^{-1}
   u_{\alpha\beta} u_{\beta\gamma} \in
   {\mathrm{Aut}}(P_\gamma)=\complex^*$ producing a \v{C}ech cocycle
   giving us the necessary data to construct a gerbe as in
   {\ref{gerbessmooth}}.
\end{proof}

\subsection{Orbifolds as Stacks}

Now we can define the stack associated to an orbifold. Let $X$ be
an orbifold with $\{(V_p,G_p,\pi_p)\}_{p \in X}$ as orbifold
structure. Let $\CC$ be the category of all open subsets of $X$
with the inclusions as morphisms and for $U \subset X$, let
$\SSS_U$ be the category of all uniformizing systems of $U$ such
that they are equivalent for every $q \in U$ to the orbifold
structure, in other words
$$\SSS_U = \left\{ (W,H,\tau) | \forall q \in U, (V_q,G_q,\pi_q) \& (W,H,\tau) \mbox{ are
equivalent at $q$} \right\}$$

By Lemma \ref{induceduniformizingsystem} and the definition of
orbifold structure it is clear that the category $\SSS$ is
fibered by groupoids. It is known, and this requires more work,
that this system $\SSS \to \CC$ is also an stack, a
Deligne-Mumford stack in the smooth case.

\subsection{Stacks as groupoids}

The following theorem has not been used in this paper, but it is
the underlying motivation for the approach that we have followed.

\begin{theorem} \cite{Fulton}  \label{stackfromgroupoid}
Every Deligne-Mumford stack comes from an \'etale groupoid scheme. Moreover, there is
a functor $\RR\rightrightarrows\UU \mapsto [ \RR \rightrightarrows \UU]$ from the
category of
 groupoids to the category of stacks inducing this realization.
\end{theorem}

When $s,t$ are smooth we can realize in a similar manner Artin
stacks \cite{Fulton}.

\bibliographystyle{amsplain}
\bibliography{gerbesbib}

\end{document}